\documentclass[12pt, leqno]{article}
\usepackage{fullpage} 
\usepackage{amsmath}
\usepackage{amssymb}
\usepackage{undertilde}
\usepackage{stmaryrd}
\usepackage{bbding}
\usepackage{enumerate}
\usepackage{graphicx}
\usepackage{xypic}
\xyoption{curve}
\usepackage{multicol}
\xyoption{color}
\usepackage{setspace}
\usepackage{boxedminipage}
\usepackage{nicefrac}
\usepackage{multicol}
\usepackage{amsthm}
  \usepackage[pdftex,
  bookmarks = false,
  pdfstartview = FitBH,
  linktocpage = true,
  pagebackref = true,
  pdfhighlight = /O,
  colorlinks=true,
  linkcolor=blue,
  citecolor=blue,
  filecolor = blue,
  urlcolor = blue,
  menucolor = blue,
]{hyperref}

\numberwithin{equation}{section}

\makeatletter
\DeclareRobustCommand\widecheck[1]{{\mathpalette\@widecheck{#1}}}
\def\@widecheck#1#2{%
    \setbox\z@\hbox{\m@th$#1#2$}%
    \setbox\tw@\hbox{\m@th$#1%
       \widehat{%
          \vrule\@width\z@\@height\ht\z@
          \vrule\@height\z@\@width\wd\z@}$}%
    \dp\tw@-\ht\z@
    \@tempdima\ht\z@ \advance\@tempdima2\ht\tw@ \divide\@tempdima\thr@@
    \setbox\tw@\hbox{%
       \raise\@tempdima\hbox{\scalebox{1}[-1]{\lower\@tempdima\box
\tw@}}}%
    {\ooalign{\box\tw@ \cr \box\z@}}}
\makeatother

\newcounter{parnum}

\newcounter{foo}

\newenvironment{myequation}
{\setcounter{equation}{\value{foo}}\begin{equation}}
{\setcounter{foo}{\value{equation}}\end{equation}\ignorespacesafterend}

\newenvironment{myeqnarray}
{\setcounter{equation}{\value{foo}}\begin{eqnarray}}
{\setcounter{foo}{\value{equation}}\end{eqnarray}\ignorespacesafterend}

\def\maxenumi{0}
\def\saveenumi{0}

\makeatletter
\AtEndDocument{%
  \immediate\write\@auxout{%
    \string\gdef\string\maxenumi{\saveenumi}}}
\makeatother

\usepackage{enumitem}

\setlist[enumerate,1]{%
  leftmargin=*,
  widest=\maxenumi,
  after=\ifnum\value{enumi}>\saveenumi\xdef\saveenumi{\the\value{enumi}}\fi}

\newenvironment{myenumerate}
{\begin{enumerate}\setcounter{enumi}{\value{foo}}}
{\setcounter{foo}{\value{enumi}}\end{enumerate}\ignorespacesafterend }

\numberwithin{equation}{section}
\numberwithin{enumi}{section}
\numberwithin{foo}{section}

\usepackage{lipsum}
\usepackage{enumitem}

\begin{document}

\title{Predicativity, the Russell-Myhill Paradox, and Church's Intensional Logic
}

\author{Sean Walsh
}

\maketitle

\begin{abstract}
This paper sets out a predicative response to the Russell-Myhill paradox of propositions within the framework of Church's intensional logic. A predicative response places restrictions on the full comprehension schema, which asserts that every formula determines a higher-order entity. In addition to motivating the restriction on the comprehension schema from intuitions about the stability of reference, this paper contains a consistency proof for the predicative response to the Russell-Myhill paradox. The models used to establish this consistency also model other axioms of Church's intensional logic that have been criticized by Parsons and Klement: this, it turns out, is due to resources which also permit an interpretation of a fragment of Gallin's intensional logic. Finally, the relation between the predicative response to the Russell-Myhill paradox of propositions and the Russell paradox of sets is discussed, and it is shown that the predicative conception of set induced by this predicative intensional logic allows one to respond to the Wehmeier problem of many non-extensions.

\end{abstract}

\tableofcontents

\section{Introduction}\label{sec:RM:01}

The Russell-Myhill paradox is a paradox about propositions which is structurally analogous  to the Russell paradox about sets. While predicativity has been well-explored as a response to the Russell paradox about sets, it seems that there has been no attempt to set out and analyze a predicative response to the Russell-Myhill paradox. The primary aim of this paper is to do just that. The crucial idea behind the predicativity response is to restrict the comprehension schema for the ambient higher-order logic. Intuitively, the comprehension schema says that every well-formed formula determines a higher-order entity. Besides the burden of showing formal consistency, a predicative response to the Russell-Myhill paradox must provide at least the beginnings of an account of why some but not all formulas succeed in determining higher-order entities.

The resulting formal system whose consistency we establish is centered around the intensional logic of Church. This intensional logic has a neutral core which we exposit in \S\ref{sec:RM:02}; it is neutral in the sense that its axioms are comparatively undemanding and consistent with contemporary theorizing based on possible-world semantics. In the subsequent \S\ref{sec:RM:03}, we set out a formalized version of the Russell-Myhill paradox, which is broadly similar to formalizations offered by Anderson and Klement, and we describe how different models offered by Kaplan and Anderson block different premises in the formalized version of the paradox. Then in \S\S\ref{sec:RM:04}-\ref{sec:RM:05}, we turn to and discuss the predicativity response to the Russell-Myhill paradox. In addition to discussing the philosophical motivations for predicative restrictions, we establish the formal consistency of the system by constructing a series of models. 

However, the models offered here validate an additional axiom of Church's intensional logic which, as Parsons and Klement have emphasized, is in some ways contrary to the spirit of a fine-grained theory of propositions. In \S\ref{sec:RM:06} we explain why this axiom holds on our models: it turns out that this is related to an expressive resource which allows these models to interpret a fragment of Gallin's intensional logic. Finally, in \S\ref{sec:RM:07}, we present one application of a broadly predicative perspective on the Russell-Myhill paradox about propositions, namely a response to the Wehmeier problem of many non-extensions that arises in connection to the naive conception of set found in Frege's {\it Grundgesetze}.

This paper is the second in a series of three papers -- the other two being \cite{Walsh2014ac}, \cite{Walsh2014ad}-- which collectively constitute a sequel to the ``Basic Law~V'' component of our earlier paper \cite{Walsh2012aa}. In the companion paper \cite{Walsh2014ac}, we use G\"odel's constructible sets to study how much of Zermelo-Fraenkel set theory can be consistently realized in these fragments of the \emph{Grundgesetze}. In the complementary paper \cite{Walsh2014ad}, we examine the deductive strength of a related theory of abstraction principles.

However, these papers do not touch the question of whether the models used to prove consistency of the {\it Grundgesetze} system are anything like intended models. In  \S\ref{sec:RM:05} of this paper, we use G\"odel's constructible sets to produce models of our extension of Church's intensional logic. Our response to the problem of many non-extensions in \S\ref{sec:RM:07} involves showing how this extension of Church's intensional logic can define a model of Frege's \emph{Grundgesetze} system. In addition to articulating a predicative response to the Russell-Myhill paradox of propositions, this paper suggests the possibility of viewing the consistent fragments of Frege's naive set theory through the lens of a consistent version of Church's intensional logic. However, despite these connections to our earlier papers, this paper has been written so that one need not have read these other papers. At one point in \S\ref{sec:RM:05} below, we reference the earlier paper \cite{Walsh2014ac} for examples of one of our defined notions in this paper (namely that of an intensional hierarchy~(\ref{eqn:RM:defn:intensional:hierarchy})). However, this is the only respect in which this paper depends on the earlier papers.

\section{The Neutral Core of Church's Intensional Logic}\label{sec:RM:02}


The intensional logic of Church is an attempt to axiomatize Frege's sense-reference distinction. Of course, Frege thought that words not only designate their referent, but also express their sense. Hence, on Frege's view, our words bear two semantic relations to non-linguistic entities, namely they bear the \emph{designation relation} to their referents and they bear the \emph{expression relation} to their senses. In the crudest of terms, Frege is a semantic dualist. This of course allowed him to say how ``the morning star'' differs in meaning from the ``evening star'': while these two linguistic expressions refer to the same planet, they express different senses.

In our view, Church's crucial observation was that semantic dualism induces a canonical non-semantic relation. The semantic dualist doesn't only think that there are more word-world connections, but they are also committed to an additional \emph{world-world} relation. In the case of Frege, the commitment is to a relationship between the abstract Fregean sense expressed by a linguistic expression and the entity (perhaps a planet) which is denoted by that linguistic expression. This relation is called the \emph{presentation} relation in the literature (\cite{Klement2010aa} p. 155), and one says that sense~$s$ \emph{presents} denotation~$d$ and one writes~$\Delta(s,d)$ precisely in the circumstance where there is a linguistic expression which expresses~$s$ and denotes~$d$. The ``triangle'' notation~$\Delta$ for the presentation relation is helpful here because it reminds us that a sense~$s$ on the bottom-left of the triangle stands in the presentation relation~$\Delta(s,d)$ to the denotation~$d$ on the bottom-right of the triangle in virtue of its semantic relations to some suppressed linguistic expression standing at the top of the triangle.\index{Symbol, Presentation~$\Delta$}

Church proceeded by axiomatizing the presentation relationship. Of course, this is not the only way that one might seek to understand the presentation relationship. With respect to a given formal language, we know how to recursively define a satisfaction relation in terms of reference, and one might have thought that one ought to proceed similarly with the presentation relationship. However, this procedure would require the notion of sense to be as conceptually transparent as the notion of reference. By contrast, Church's own aim in axiomatizing the presentation relationship was to dissipate outright skepticism about Fregean sense. Here is how Church put the point in a 1943 review of a paper of Quine:
\begin{quote}
There remains the important task, which has never been approached, of constructing a formalized semantical system which shall take account of both kinds of meaning, the relation between a name and its denotation, and the relation between a name and its sense. [\ldots] [\P] [\ldots] Ultimately it is only on the basis of their inclusion in an adequate system of this kind that such otherwise indefensibly vague ideas as `understanding' of an expression, `attribute,' `objectiver Inhalt des Denkens,' may be regarded as logically significant (\cite{Church1943aa} p. 47).
\end{quote}
Hence one of the original aims of Church's work was to produce a formal theory whose quantifiers ranged over Fregean senses and which thus serve to implicitly define the notion of Fregean sense.

The first component of Church's formal theory concerned the Fregean doctrine that sense determines reference. Frege tells us that whenever two linguistic expressions have the same sense, then if one refers then the other does too and they have the same referent. When put this way, it automatically suggests the following axiom (Church's Axiom~17, cf. \cite{Church1951ab} p.~19, \cite{Klement2002aa} pp.~108-109, \cite{Anderson1980aa} p. 220; \cite{Anderson1984aa} Axiom C8 p. 377):
\begin{myenumerate}
\item \emph{Sense Determines Reference}:~$(\Delta(s,d_0) \; \& \; \Delta(s,d_1)) \Longrightarrow d_0=d_1$ \label{eqn:RM:sdr}
\end{myenumerate}
Practically, this indicates to us that the presentation relationship is functional in character.  Thus instead of writing~$\Delta(s,d)$, we may write instead~$\Delta(s)=d$. Likewise, borrowing notation from computability theory, sometimes we write~$\Delta(s)\hspace{-1mm}\downarrow$ to indicate that there is a~$d$ such that~$\Delta(s)=d$ (cf. \cite{Soare1987} pp. 16-17). Of course, we should keep in mind that on its intended interpretation, the presentation relation is not a total function. For, any meaningful linguistic expression will always have a sense but need not have a referent.

The other of Church's axioms pertain to compositionality. On the side of reference, Frege postulated a fundamental distinction between objects and concepts. Objects were the referents of proper names of people and places, whereas concepts were the referents of predicate-words. In sentences such as ``Venus is a planet,'' we predicate a concept (``being a planet'') of an object (``Venus''). This may be rendered as a case of functional application, namely~$\textsc{Planet}(\textsc{Venus})=1$. Following contemporary practice, we here identify ``$1$'' with the truth-value ``true'' and ``$0$'' with the truth-value ``false,'' and for the sake of simplicity we assume that these are the only truth-values.

Due to the fact that sense is a determiner of reference, it's natural to think that senses of sentences are also compositional. Frege called the senses of sentences \emph{thoughts}, so that sentences express thoughts and refer to truth-values. Just as the reference of a sentence is a function of the reference of its constituent parts, so Frege and Church hold that the sense of a sentence (a thought or a proposition) is a function of the senses of its parts. This thus suggested to Church the following axiom on the presentation relationship (Church's Axiom 15 \cite{Church1951ab} p.~18, \cite{Klement2002aa} pp.~108-109, \cite{Anderson1980aa} p. 219; cf. \cite{Anderson1984aa} Axiom C7 p. 377):
\begin{myenumerate}
\item  \emph{Composition Axiom}:~$[\Delta(f^{\prime})=f \; \& \; \Delta(x^{\prime})=x] \Longrightarrow \Delta(f^{\prime}\langle x^{\prime}\rangle) = f(x)$\label{eqn:RM:cca1}
 \end{myenumerate}
 Here~$(f^{\prime}, x^{\prime})\mapsto f^{\prime}\langle x^{\prime}\rangle$ is a primitive \emph{intensional} application function on senses,\index{Symbol, Intensional Application~$(f^{\prime}, x^{\prime})\mapsto f^{\prime}\langle x^{\prime}\rangle$} just as~$(f,x)\mapsto f(x)$ is a primitive \emph{extensional} application function on referents.\index{Symbol, Extensional Application~$(f,x)\mapsto f(x)$} The axiom itself leaves open the relationship between intensional and extensional application, although we'll see later in this section, Church himself proposed that we identify them. 

Now we have at least six types of entities: objects, senses of objects, concepts, senses of concepts, truth-values, senses of truth-values (which we also call \emph{thoughts} or \emph{propositions}). However, there is a serious redundancy here. For, we may identify concepts with functions from objects to truth-values. (For details on this familiar identification, see circa equation~(\ref{eqn:RM:ide}) in the next section). If one does so, then it seems natural enough to further assume that if~$a$ is a type of entity and~$b$ is a type of entity, then there is a type~$ab$ of entities consisting of functions from entities of type~$a$ to entities of type~$b$. One makes these assumptions rigorous by defining the types recursively as follows:  
\begin{myenumerate}
\item (Types in the Church System) (i) there is a type~$e$ of objects, (ii) there is a type~$t$ of truth-values, (iii) if~$a,b$ are types, then there is a type~$ab$ of functions from type~$a$ entities to type~$b$ entities, and (iv)~if~$a$ is a type then~$a^{\prime}$ is a type of senses which present entities of type~$a$.\label{eqn:typesystemChurch}\index{Types in Church's System~(\ref{eqn:typesystemChurch})}
\end{myenumerate}
\noindent In this last clause, it's important to emphasize that~$a\mapsto a^{\prime}$ is a primitive operation on types (cf. Kaplan~\cite{Kaplan1975aa} p. 721 and Klement~\cite{Klement2010aa} p. 173). Hence,~$a^{\prime}$ is the result of applying an operation to the type~$a$, and not simply another variable for types (and likewise~$a^{\prime\prime}$ is the result of applying the prime operation to type~$a^{\prime}$). Sometimes in what follows, if we write entities of type~$a$ as~$f,g,h,\ldots$ (resp.~$x,y,z \ldots$), then we will adopt the convention of writing entities of type~$a^{\prime}$ as ~$f^{\prime},g^{\prime},h^{\prime},\ldots$ (resp.~$x^{\prime},y^{\prime},z^{\prime} \ldots$). However, under this convention, entity~$f^{\prime}$ of type~$a^{\prime}$ is \emph{not} the result of applying any operation to the entity~$f$ of type~$a$, but rather just a conventional device which allows us to visually keep track of which entity has which type. 

Having set up the type system in this way, one sees immediately that there must be not a single presentation relationship~$\Delta$, but rather a presentation relationship~$\Delta_a$\index{Symbol, Presentation Symbol, Typed~$\Delta_a$} for each type~$a$, which relates senses of type~$a^{\prime}$ to entities of type~$a$. Having made this distinction, one thus reformulates the Axiom that Sense Determines Reference~(\ref{eqn:RM:sdr}) and the Composition Axiom~(\ref{eqn:RM:cca2}) as follows:
\begin{myenumerate}
\item \emph{Typed Sense Determines Reference}:~$(\Delta_a(s,d_0) \; \& \; \Delta_a(s,d_1)) \Rightarrow d_0=d_1$ \label{eqn:RM:sdr2}\index{Sense Determines Reference, Typed~(\ref{eqn:RM:sdr2})}
\item  \emph{Typed Composition}:~$[\Delta_{ab}(f^{\prime})=f \; \& \; \Delta_a(x^{\prime})=x] \Longrightarrow \Delta_b(f^{\prime}\langle x^{\prime}\rangle) = f(x)$\label{eqn:RM:cca2}\index{Composition, Typed~(\ref{eqn:RM:cca2})}
 \end{myenumerate}
\noindent In the latter, the intensional application function~$(f^{\prime}, x^{\prime})\mapsto f^{\prime}\langle x^{\prime}\rangle$ takes an sense~$f^{\prime}$ of type~$(ab)^{\prime}$ and a sense~$x^{\prime}$ of type~$a^{\prime}$ and returns a sense of type~$b^{\prime}$, just as the extensional application function~$(f,x)\mapsto f(x)$ takes a referent~$f$ of type~$ab$ and a referent~$x$ of type~$a$ and returns a referent of type~$b$. Hence, just as there are as many presentation relations as there are types, so there are as many intensional and extensional application functions as there are pairs of types. Later, when we deal more formally with these systems, we will introduce symbols subscripted by types for the intensional and extensional application functions but for the time being we simply allow context to determine the types (cf. circa equation~(\ref{eqn:RM:defn:intensional}) in \S\ref{sec:RM:05}). Finally, let's record our standing assumption in this paper that entities of type~$ab$ are individuated extensionally:
\begin{myenumerate}
\item \emph{Extensional Identity Criterion for Functional Entities}: if~$f,g$ are entities of type~$ab$ then~$f=g$ if and only if~$f(x)=g(x)$ for all entities~$x$ of type~$a$.\label{extensionalidentityfunctions}\index{Extensional Identity for Functional Entities}
\end{myenumerate}
But no analogous assumptions are made on the individuation of entities of type~$(ab)^{\prime}$ in this paper.

In his own writings, Church adopted the following axiom on the types themselves:
\begin{myenumerate}
\item \emph{Axiom of Type Reduction}:~$(ab)^{\prime} = a^{\prime}b^{\prime}$\label{eqn:RM:TRD}\index{Axiom of Type Reduction~(\ref{eqn:RM:TRD})}
\end{myenumerate}
In other words, this axiom says that the type~$(ab)^{\prime}$ of senses of functions from type~$a$ entities to type~$b$ entities is identical to the type~$a^{\prime}b^{\prime}$ of functions from senses of type~$a$ entities to senses of type~$b$ entities. It's called a reduction axiom because it allows one to reduce the senses of all higher-order entities to the senses of objects and truth-values (and senses of senses of objects, senses of propositions, etc.). The primary formal advantage of doing this is that it allows one to interpret intensional application as extensional application and hence disburdens one from developing an alternative conception of intensional application. Indeed, if senses of concepts are really extensional functions from senses of objects to senses of truth-values, then it's natural to think that the sense of a proposition is produced via the extensional application of the sense of an object to the sense of a concept.

But some of the objections to Church's intensional logic have revolved around this Axiom of Type Reduction~(\ref{eqn:RM:TRD}). Dummett was concerned that the reduction axiom would require us to deny the seemingly plausible idea that ``[\ldots] we are able to learn what thought \emph{some} sentences containing the predicate express in advance of knowing the sense of the predicate'' (\cite{Dummett1981aa} p. 294, cf. \cite{Klement2002aa} pp. 69-70 ff). Dummett's idea was that we can learn the senses of complete propositions like~$Fa$ and~$Fb$ without precisely knowing the sense of~$F$. But the Axiom of Type Reduction~(\ref{eqn:RM:TRD}) demands that the sense of~$F$ is a function, hence presumably complete propositions like~$Fa$ and~$Fb$ will be the result of functional application of this sense-function, so that knowledge of them may well require prior knowledge of the sense of~$F$.

Another objection to Church's Axiom of Type Reduction~(\ref{eqn:RM:TRD}) is due to Bealer, who notes that functionality seems to be entirely absent from qualia and other facets of conscious experience. Bealer writes: ``Joy, the shape of my hand, the aroma of coffee- these are not functions. When I feel joy, see the shape of my hand, or smell the aroma of coffee, it is not a function that I feel, see, or smell (cf. \cite{Bealer1982aa} p. 90, cf. \cite{Duzi2010aa} p. 6). This concern resonates well with the observation that when we are pressed to say something about the sense of words such as ``red,'' ``cold,'' or ``bitter,'' the mathematical notion of a function is far from our first thoughts. Of course, someone who denies Church's reduction axiom for these types of reasons need not be taken to deny that senses can compose with other senses. Rather, the denial should be registered merely as a denial that senses of predicate words can be exhaustively identified with functions from senses to other senses.

Once we reject Church's reduction axiom, we are left with the following core of Church's system: 
\begin{myenumerate}
\item The \emph{core of Church's system} consists of the Typed Sense Determines Reference~(\ref{eqn:RM:sdr2}) and the Typed Composition Axiom~(\ref{eqn:RM:cca2}); this theory is a typed theory, and the types are exactly as in~(\ref{eqn:typesystemChurch}). \label{eqn:RM:coresystem}\index{Core of Church's System~(\ref{eqn:RM:coresystem})}
\end{myenumerate}
\noindent It's noteworthy that there is nothing in these core axioms themselves that forces or even necessarily recommends the identification of type~$a^{\prime}$ with Fregean senses as opposed to any other notion of meaning. Indeed, Kaplan pointed out long ago that the standard frameworks for possible worlds semantics yield models of these axioms (cf. \cite{Kaplan1975aa} pp. 721~ff). In particular, Kaplan proceeded by identifying type~$a^{\prime}$ in the Church system with the type of functions from the worlds to the entities of type~$a$ and by defining the other primitives of the Church system as follows, wherein~$w_0$ is fixed world (say the actual world) and~$w$ is an arbitrary world:
\begin{myequation}\label{eqn:kapalandadsfas}
\Delta_a(f^{\prime}) = f^{\prime}(w_0), \hspace{10mm} (f^{\prime}\langle x^{\prime}\rangle)(w) =  (f^{\prime}(w))(x^{\prime}(w))
\end{myequation}
On the basis of these definitions, it's not too difficult to check that the Typed Sense Determines Reference~(\ref{eqn:RM:sdr2}) and the Typed Composition Axiom~(\ref{eqn:RM:cca2}) are both satisfied.\footnote{The Typed Sense Determines Reference~(\ref{eqn:RM:sdr2}) is satisfied because the definition of~$\Delta_a(f^{\prime})$ in equation~(\ref{eqn:kapalandadsfas}) is clearly functional since~$f^{\prime}$ is by stipulation a function defined on worlds and the world of evaluation~$w_0$ is fixed. For the Typed Composition Axiom~(\ref{eqn:RM:cca2}), suppose that~$\Delta_{ab}(f^{\prime})=f$ and~$\Delta_a x^{\prime}=x$. Then one can calculate that~$\Delta_b(f^{\prime}\langle x^{\prime}\rangle) = (f^{\prime}\langle x^{\prime}\rangle)(w_0) = (f^{\prime}(w_0))(x^{\prime}(w_0))=(\Delta_{ab} f^{\prime})(\Delta_a x^{\prime}) = f(x)$.} The core of Church's system is thus fairly neutral on the philosophical interpretation of the intensional notions which it axiomatizes. In the next section we turn to the formalization of the Russell-Myhill paradox within this core~system, and this treatment of the paradox will yield a common framework in which advocates of Fregean sense and advocates of possible worlds semantics can discuss the comparative advantages and disadvantages of different solutions to the paradox.

Before doing so, it's perhaps worth underscoring some of the departures that we have made in this paper from traditional treatments of Church's intensional logic. Some of these differences are merely notational. One such difference is that Church wrote the type reserved for functions from entities of type~$a$ to entities of type~$b$ as~$ba$ rather than~$ab$ (cf. \cite{Church1951ab} p. 12, \cite{Anderson1984aa} p. 370, \cite{Klement2002aa} p. 106). We prefer the latter simply because it is now the norm in formal semantics (cf. \cite{Heim1998} p. 28, \cite{Gamut1991aa} pp. 84, 121). Further, Church respectively used the letters~$o_1$ and~$\iota_1$ instead of~$t$ and~$e$ for the truth-values and objects (cf. \cite{Church1951ab} p. 11, \cite{Anderson1984aa} p. 370, \cite{Klement2002aa} p. 106). Again, we use~$t$ and~$e$ simply because this is now the norm (cf. \cite{Heim1998} p. 28, \cite{Gamut1991aa} pp. 79, 128). Finally, sometimes in Church --and sometimes in intensional logics more generally-- the word ``concept,'' perhaps proceeded by modifiers like ``individual'' or ``propositional,'' is reserved for certain kinds of senses or intensions (cf. \cite{Klement2002aa} p. 96, \cite{Duzi2010aa} pp. 155 ff, \cite{Gamut1991aa} p. 122). However, here in this paper we eschew this usage and use the Fregean terminology, on which ``concepts'' are the unsaturated entities which may be saturated by objects and which are thus one-half of the concept-object distinction. 

The chief contentful difference between Church's original formulation of his intensional logic and our treatment of it concerns intensional application. Church himself did not introduce a primitive~$(f^{\prime},x^{\prime})\mapsto f^{\prime}\langle x^{\prime}\rangle$ for intensional application; again, this was because of his adoption of the axiom on type-reduction~(\ref{eqn:RM:TRD}) which we and many others reject. Further, because of his adoption of this axiom, Church did not introduce a primitive operation~$a\mapsto a^{\prime}$ on types but could simply get by with postulating types~$o_1, o_2, \ldots, \iota_1, \iota_2, \ldots$ wherein the type $o_1$ is the type of truth-values and the type $\iota_1$ is the type of objects and wherein~$\sigma_{n+1}$ is the type reserved for senses of entities of type~$\sigma_n$ for~$\sigma\in \{o,\iota\}$ (cf. \cite{Church1951ab} pp. 7, 11, \cite{Anderson1984aa} p. 370, \cite{Klement2002aa} p. 106). If one rejects the axiom on type-reduction~(\ref{eqn:RM:TRD}), then it's natural to postulate the primitive operation~$a\mapsto a^{\prime}$ on types, and here we follow Kaplan~\cite{Kaplan1975aa} p. 721 and Klement~\cite{Klement2010aa} p. 173. Finally, much of Church's own work on his system concerned various proposals for individuating senses, and these went under the name of Alternative~(0), Alternative~(1), and Alternative~(2) (cf. Klement~\cite{Klement2002aa} pp. 101~ff for overview). How exactly senses are individuated is obviously important, but it is not needed for the formalization of the Russell-Myhill Paradox discussed in the next section or for the predicative response discussed in the subsequent sections.\footnote{Admittedly, there is something deeper going on here. The distinction between Alternative~(0) and Alternative~(1) lies in whether lambda-conversion preserves sense. However, lambda-terms are an alternative way of formalizing comprehension~(cf. (\ref{eqn:RM:compschem})) which the predicative response offered here does not have available in full generality. Thus the formal extensions of Church's core system with which we work here simply don't have lambda-terms in the object-language. Hence an immediate issue which faced Church-- namely whether to say that lambda-conversion preserves sense-- is not even available in the object-language of our systems.}
  
\section{Formalized Version of the Russell-Myhill Paradox and Extant Responses}\label{sec:RM:03}

The aim of this section is to set out a formalized version of the Russell-Myhill Paradox and to survey some extant non-predicative responses. The formalization offered here is distinct from but owes much to the formalizations offered by Anderson and Klement, and we'll discuss explicitly in this section these similarities and differences. Since the Russell-Myhill paradox is a proposition-theoretic version of the Russell paradox about sets, it's useful to begin with a brief review of Russell's paradox about sets and Cantor's related theorem about cardinalities. A collection~$X$ is said to have \emph{cardinality less than or equal} to collection~$Y$ just in case there is an injection~$\iota:X\rightarrow Y$, while the two collections~$X,Y$ are said to have \emph{the same cardinality} if there is a bijection between them. Cantor's theorem about cardinalities says that for any collection, there is no injection from the set of all its subcollections to itself. In symbols, Cantor's theorem says that for any~$X$, there is no injection from~$\{Y: Y\subseteq X\}$ to~$X$ itself. But there is a natural bijective correspondence between the subcollections~$Y\subseteq X$ and the zero-one valued functions~$f:X\rightarrow \{0,1\}$, given by sending~$Y\subseteq X$ to its characteristic function~$f_Y:X\rightarrow \{0,1\}$ which is defined by
\begin{myequation}\label{eqn:RM:ide}
f_Y(x) = \begin{cases}
1      & \text{$x\in Y$}, \\
0     & \text{$x\notin Y$}.
\end{cases}
\end{myequation}
Since the type~$at$ is reserved for functions from entities of type~$a$ to the truth-values~$\{0,1\}$, there is thus a natural type-theoretic expression of Cantor's theorem:
\begin{myenumerate}
\item (Type-Theoretic Version of Cantor's Theorem) For any type~$a$, there is no injection from entities of type~$at$ to entities of type~$a$.\label{eqn:RM:typetheoreticRP}\index{Type-Theoretic Version of Cantor's Theorem~(\ref{eqn:RM:typetheoreticRP})}
\end{myenumerate}
This version is entirely type-theoretic, since the injection in question would be an element of type~$(at)a$ and since the property of being injective is expressible purely with the extensional application notions built into the type theory.

Let us briefly recall the traditional proof of the type-theoretic version of Cantor's Theorem. For the sake of readability, in this proof let us write entities of type~$at$ as~$f,g,h$ \ldots and entities of type~$a$ as~$x,y,z\ldots$. A function from the entities of type~$at$ to entities of type~$a$ is a function~$\iota$ taking input~$f$ of type~$at$ and returning output~$\iota(f)=x$ of type~$a$. Now, suppose that there was such an injection~$\iota$ from entities of type~$at$ to entities of type~$a$. Then consider the \emph{diagonal} map~$d$ from elements of type~$a$ to elements of type~$t$ given by\index{Diagonal map~$d$ (\ref{eqn:RM:eqndai})}
\begin{myequation}\label{eqn:RM:eqndai}
d(x) = 1 \Longleftrightarrow \exists \; f \; (\iota(f) = x \; \& \; f(x)=0)
\end{myequation}
Then let~$y=\iota(d)$ and ask whether~$d(y)=1$ or~$d(y)=0$. If~$d(y)=1$, then by the left-to-right direction of equation~(\ref{eqn:RM:eqndai}) one obtains witness~$f$ satisfying~$\iota(f) =y$ and~$f(y)=0$. Then since~$\iota(d)=y=\iota(f)$, we may conclude from the injectivity of~$\iota$ that~$d=f$, which contradicts that~$d(y)=1$ while~$f(y)=0$. Alternatively, if~$d(y)=0$, then~$d$ and~$y$ are witness to the right-hand side of equation~(\ref{eqn:RM:eqndai}), and so by the right-to-left direction of this equation we have~$d(y)=1$, a contradiction. Hence, in either case we obtain a contradiction. 

The connection between the type-theoretic version of Cantor's theorem and Russell's paradox about sets can be made more transparent if one defines a membership relation by
\begin{myequation}
y\in x \Longleftrightarrow \exists \; f \; (\iota(f)=x \; \& \; f(y)=1)
\end{myequation}
Further, for the moment let's call a \emph{set} an entity which is in the range of the operator~$\iota$. Then for any set~$x$, it follows from the injectivity of the~$\iota$-operator that we can also express \emph{non}-membership in~$x$ with an existential quantifier as follows:
\begin{myequation}
y\notin x \Longleftrightarrow \exists \; f \ (\iota(f)=x \; \& \; f(y)=0)
\end{myequation}
But then it is easy to see that for sets~$x$, equation~(\ref{eqn:RM:eqndai}) is equivalent to:
\begin{myequation}\label{eqn:RM:eqndai2}
d(x) = 1 \Longleftrightarrow x\notin x
\end{myequation}
Expressed in these terms, the diagonal function~$d$ from the above paragraph is the characteristic function of the collection of sets which are not members of themselves. This is one way to see the connection between the Russell paradox about sets and the type-theoretic version of Cantor's Theorem.

The Russell-Myhill paradox about propositions proceeds by arguing, based on considerations related directly to propositions, that there is an injection from collections of propositions to propositions. Since type~$t^{\prime}$ is reserved for propositions in the Church system and since collections of propositions can be identified with their characteristic functions of type~$t^{\prime}t$ (\emph{\`a~la} equation~(\ref{eqn:RM:ide})), if this argument were correct than it would mean that there was an injection from entities of type~$t^{\prime} t$ to entities of type~$t^{\prime}$, which contradicts the type-theoretic version of Cantor's Theorem~(\ref{eqn:RM:typetheoreticRP}).  One way to respond to the Russell-Myhill paradox about propositions is to block \emph{in a well-motivated way} the type-theoretic version of Cantor's Theorem~(\ref{eqn:RM:typetheoreticRP}). This is the kind of proposal which we shall pursue in this paper, beginning in \S\ref{sec:RM:04}. But in the remainder of this section we focus simply on formalizing the Russell-Myhill paradox and on surveying other extant responses.

The informal version of the Russell-Myhill paradox was initially described in Appendix~B of Russell's 1903 \emph{The Principles of Mathematics} and then again in Myhill's~1958 paper on Church's intensional logic.\footnote{More specifically see \S{500} p.~538 of Russell \cite{Russell1903aa} and p. 82 of Myhill \cite{Myhill1958aa}. According to the history as set out in de~Rouilhan \cite{Rouilhan2005aa}, Russell never mentioned this paradox again. As for Myhill, in the same 1958 paper he reports that Carnap's ``general approach to the problem, in terms of `possible worlds' and state-descriptions, is in [his] opinion practically certain to yield a correct explication within a few years'' (\cite{Myhill1958aa} p. 81). This contrasts with Myhill's earlier 1952 paper on Church (\cite{Myhill1952ae}) in which he weighs carefully the costs and benefits of Fregean and modal approaches without indicating a decisive preference for either. It is well-known that Myhill continued to work on intuitionistic and non-classical approaches to the set-theoretic paradoxes throughout his career, but to my knowledge he never after the 1958 paper returned to this proposition-theoretic paradox.} The argument of Russell and Myhill runs as follows. Given a collection of propositions~$\mathcal{C}$, consider the proposition~$\iota(C)$ expressed by the sentence ``every proposition in~$\mathcal{C}$ is true'' (or ``every proposition is in~$\mathcal{C}$.'') It seems that this function is an injection. For, suppose that~$\iota(\mathcal{C})=\iota(\mathcal{D})$. Since these two propositions differ only as to~$\mathcal{C}$ and~$\mathcal{D}$, then presumably the senses (or intensions) of their constituents~$\mathcal{C}$ and~$\mathcal{D}$ are the same as well. And this would presumably imply that~$\mathcal{C}$ and~$\mathcal{D}$ are the same not only in sense or intension, but that they are also the same in reference or extension. Hence the map~$\iota$ is ostensibly an injection from collections of propositions to propositions. But, by applying the type-theoretic version of Cantor's Theorem~(\ref{eqn:RM:typetheoreticRP}), we obtain a contradiction.

While this version of the Russell-Myhill paradox is traditional, it is not obviously a formal paradox. This is for two reasons. First, formal paradoxes show that some formal system is inconsistent. But it is not at all obvious-- based merely on its informal description-- in what system the argumentation of the above paragraph may be formalized. To be sure, a good start could be made on this to the extent that one could formalize the notion of ``a sentence expressing a proposition.'' But to the extent that one could formalize this notion one could presumably likewise formalize the notion of ``a sentence denoting a truth-value,'' and hence one would worry that this formalization would require prior treatment of the liar paradox.\footnote{Intensional logics like Church's intensional logic and possible worlds semantics have resources for axiomatizing the notion of a ``proposition denoting the true.'' In Church's system, this is written as~$\Delta_t(p)=1$ while in possible worlds semantics this is written~$p(w_0)=1$ where~$w_0$ is the world of evaluation. However, in neither of these intensional logics does one have the resources for going from a name of a sentence to the proposition expressed by the sentence. If one did, then since these systems of intensional logic are consistent with the addition of resources needed to effect self-reference, one could replicate the formal versions of the liar paradox.} Second, formal paradoxes are always \emph{valid} arguments, whose conclusion is that some formal set of axioms is inconsistent. But there is a real concern about the validity of the above rendition of the Russell-Myhill paradox. For, it seems at crucial points to equivocate between the collection~$\mathcal{C}$ and a sense thereof. Indeed, it seems that it is the latter which would contribute to the proposition expressed by ``every proposition in~$\mathcal{C}$ is true.'' Yet, the argument as a whole pertains to an injection which takes as inputs collections of propositions~$\mathcal{C}$. A truly formalized version of the Russell-Myhill paradox would leave no doubt as to whether the argument was, at any juncture, operating on a collection of propositions or a sense.\footnote{Klement suggests that this kind of concern is one way of understanding Frege's own reservations about the Russell-Myhill paradox (\cite{Klement2002aa} p. 183).}

The formalization of the Russell-Myhill paradox which we adopt avoids these two problems, and reads as follows:
\begin{myenumerate}
\item (\emph{Formalized Russell-Myhill Paradox}). The following axioms are jointly inconsistent against the background of the core of Church's system~(\ref{eqn:RM:coresystem}): the Surjectivity Axiom~(\ref{eqn:RM:SM}), the Senses are Objects Axiom~(\ref{eqn:RM:SO}), the Propositions as Fine-Grained as Objects Axiom~(\ref{eqn:RM:derepropsaxiom}), and the type-theoretic version of Cantor's Theorem~(\ref{eqn:RM:typetheoreticRP}).\label{eqn:RM:formalized}\index{Formalized Russell-Myhill Paradox (\ref{eqn:RM:formalized})}
\end{myenumerate}
As one can see, this formalization concerns three additional axioms which we need to introduce and motivate in these next pages prior to setting out the derivation of the paradox. 

The first axiom in the formalized version of the paradox is the called the Surjectivity Axiom. In essence, this axiom says that every entity-- including the higher-order ones-- are presented by some sense or intension:
\begin{myenumerate}
\item \emph{Surjectivity Axiom}: for each type~$a$ and each element~$f$ of type~$a$, there is element~$f^{\prime}$ of type~$a^{\prime}$ such that~$\Delta_a(f^{\prime})=f$.\label{eqn:RM:SM}\index{Surjectivity Axiom (\ref{eqn:RM:SM})}
\end{myenumerate}
The immediate warrant for this axiom is that there is simply no other way to formalize the Russell-Myhill paradox. For, consider again how it opens: ``for each collection of propositions~$\mathcal{C}$, consider the proposition~$\iota(C)$ expressed by the sentence `every proposition in~$\mathcal{C}$ is true.'~'' We accordingly need some way to move from \emph{any} collection of propositions to a proposition. It seems that any way in which we do this will take a collection of propositions, take a sense or intension which presents this collection, and build a proposition based off of this sense.\footnote{The Surjectivity Axiom has a long and complicated history in Church's own writings. In 1946, Church seemed to indicate that Cantor-like paradoxes would lead one to deny this axiom (\cite{Church1946aa} p. 31). In 1974, Church indicated that this axiom followed from the premises of his system called Alternative~2 (\cite{Church1974aa} p. 145). In his last paper in 1993, Church included this axiom in his system (\cite{Church1993aa} pp. 144-145), albeit without saying anything explicit about his reasons for this inclusion. For other statements of the Surjectivity Axiom in the secondary literature, see Anderson~\cite{Anderson1980aa} principle~(C) p. 221 and Klement \cite{Klement2002aa} Theorem LSD(0) 1 p. 116 and Klement \cite{Klement2003aa} p. 305 Axiom~PC. For more on Anderson and Klement on the Surjectivity Axiom, see the discussion below.
}

The next axiom concerns the location of senses or intensions within Frege's concept-object distinction, or within the typed systems usually employed in formal semantics.  In essence, it says that senses or intensions fall on the object side of the concept-object distinction:
\begin{myenumerate}
 \item \emph{Senses are Objects Axiom}: for each type~$a$ and each element~$f^{\prime}$ of type~$a^{\prime}$, there is element~$x$ of type~$e$ such that~$f^{\prime}=x$.\label{eqn:RM:SO}\index{Senses are Objects Axiom (\ref{eqn:RM:SO})}
\end{myenumerate}
This axiom is non-trivial because~$f^{\prime}$ and~$x$ are variables of different types. If contrary to fact they were variables of the same type, this would simply be a truth of the ambient predicate logic. The primary kind of consideration which points in favor of the Senses are Objects Axiom is a reflection on traditional conceptions of the nature of Fregean senses: in particular, while Russell suggested that we might view Fregean senses as definite descriptions,\footnote{It is admittedly somewhat inaccurate to speak of definite descriptions merely as an ``interpretation of Fregean sense,'' since they in fact provide a systematic way of dispensing with the Fregean notion of expression altogether and maintaining that reference is the sole semantic primitive. But presumably part of our tradition's reason for thinking that Frege's theory of meaning is susceptible to modal counterexamples couched in terms of definite descriptions is something like the thought that we can think of Fregean senses as definite descriptions.} Dummett has suggested that we might understand them as certain kinds of procedures or algorithms, a ``route to reference.''\footnote{Cf. \cite{Dummett1981aa} pp. 96, 102, 179~ff, \cite{Horty2007aa} pp. 66~ff, \cite{Taschek2010aa} p. 323. This idea is also associated to Tich\'y. See in particular the papers ``Sense and Procedure'' and ``Intensions in Terms of Turing Machines'' in \cite{Tichy2004aa}.} If either of these two traditional proposals about the nature of Fregean sense are correct, then it seems that senses might be regarded as objects of certain kinds, as opposed to concepts: for, whatever the exact nature of definite descriptions and algorithms, presumably they are unsaturated and fall on the ``object'' side of Frege's concept-object distinction. For instance, if one views definite descriptions as G\"odel numbers of formulas or if one views algorithms as indexes of Turing machines this will be the case.\footnote{Obviously, one way to respond to the version of the Russell-Myhill paradox formalized here would be to deny the Senses are Objects Axiom (\ref{eqn:RM:SO}). One way to do that might be to accept that senses are definite descriptions or procedures but to deny that these can be identified with specific objects like G\"odel numbers of formulas or indexes of Turing machines. Traditional reasons for such a denial might be that e.g. abstract procedures aren't represented by a specific index for a specific Turing machine, but rather by a large class of such indexes (cf. \cite{Blass2009ab}). I don't think that such a response would ultimately succeed. For, grant all this and then just select, for each abstract procedure, a specific index for a specific Turing machine which represents it, and call these things \emph{quasi-senses}. Then quasi-senses are objects and so one could run the entire Russell-Myhill paradox again with respect to quasi-senses. For, the other axioms occurring in the formalized version of the paradox seem just as plausible for the so-defined quasi-senses as for senses qua abstract procedures. A similar point can be made with respect to definite descriptions simply by selecting G\"odel numbers of specific formulas.}

The final axiom operative in our formalized version of the Russell-Myhill paradox~(\ref{eqn:RM:derepropsaxiom}) is an axiom postulating a connection between objects and propositions:
\begin{myenumerate}
\item \emph{Propositions as Fine-Grained as Objects Axiom}: there is an injection~$\chi$ from entities of type~$e$ to entities of type~$t^{\prime}$. \label{eqn:RM:derepropsaxiom}\index{Propositions as Fine-Grained as Objects Axiom (\ref{eqn:RM:derepropsaxiom})}
\end{myenumerate}
Since entities of type~$e$ are objects and entities of type~$t^{\prime}$ are propositions, this axiom is just saying that there is an injection from objects to propositions. One plausible case for this axiom might be made from the assumption that (i)~our language is ample enough to distinguish different objects from one another and (ii)~propositions are organized roughly after the manner of the sentences which express them. For, by (i), for any object, we can fasten onto a predicate or name in our language which distinguishes this object from the others in our purview. And then by (ii) the distinctness of this item of language, be it a predicate or name, will then be matched in the propositions expressed by sentences featuring it.

Having set out and motivated the various axioms, let us now establish the formalized version of the Russell-Myhill paradox~(\ref{eqn:RM:formalized}). By the Surjectivity Axiom~(\ref{eqn:RM:SM}), for every collection of propositions~$\mathcal{C}$, there is a sense~$\mathcal{C}^{\prime}$ which presents it. By the Senses are Objects Axiom~(\ref{eqn:RM:SO}), each such sense~$\mathcal{C}^{\prime}$ is identical to some object. There is thus a map from collections of propositions to objects such that the object is identical to a sense which presents the collection. More formally:
\begin{myenumerate}
\item For every collection of propositions~$\mathcal{C}$ there is an object~$x$ and there is a sense~$\mathcal{C}^{\prime}$ such that~$\mathcal{C}^{\prime}$ is identical to object~$x$ and~$\mathcal{C}^{\prime}$ presents~$\mathcal{C}$.\label{eqn:RM:AC1}
\end{myenumerate}
This induces a map~$\mathcal{C}\mapsto \delta(\mathcal{C})$ from collections of propositions to objects such that 
\begin{myenumerate}
\item For every collection of propositions~$\mathcal{C}$ there is a sense~$\mathcal{C}^{\prime}$ such that~$\mathcal{C}^{\prime}$ is identical to object~$\delta(\mathcal{C})$ and~$\mathcal{C}^{\prime}$ presents~$\mathcal{C}$.\label{eqn:RM:AC2}
\end{myenumerate}
Further, by the Propositions as Fine-Grained as Objects Axiom~(\ref{eqn:RM:derepropsaxiom}), there is an injection~$\chi$ from objects to propositions. Let~$\iota$ be the composition of the two maps, so that~$\iota = \chi \circ \delta$. Then~$\iota$ is a map from collections of propositions~$\mathcal{C}$ to propositions~$\iota(\mathcal{C})$. Further, the map~$\iota$ too is an injection. For, suppose that~$\iota(\mathcal{C}_1)=\iota(\mathcal{C}_2)$. Then since~$\chi$ is an injection,~$\delta(\mathcal{C}_1)=\delta(\mathcal{C}_2)$. Then by the characterization of~$\delta$ in equation~(\ref{eqn:RM:AC2}), for each~$k\in\{1,2\}$ there is sense~$\mathcal{C}^{\prime}_k=\delta(\mathcal{C}_k)$ which presents~$\mathcal{C}_k$. But since~$\delta(\mathcal{C}_1)=\delta(\mathcal{C}_2)$, we have that~$\mathcal{C}_1^{\prime}=\mathcal{C}_2^{\prime}$. Since~$\mathcal{C}_k^{\prime}$ presents~$\mathcal{C}_k$, by the Typed Sense Determines Reference Axiom~(\ref{eqn:RM:sdr2}), it follows that~$\mathcal{C}_1=\mathcal{C}_2$. Hence the map~$\iota$ is an injection from collections of propositions to propositions, which contradicts the type-theoretic version of Cantor's Theorem~(\ref{eqn:RM:typetheoreticRP}).

This formalization of the Russell-Myhill paradox is distinct from but owes much to the earlier formalizations of Anderson and Klement. On the one hand, Anderson made the Surjectivity Axiom~(\ref{eqn:RM:SM}) the focus of his treatments of the paradox (\cite{Anderson1980aa} pp. 221~ff, \cite{Anderson1987aa} pp. 107~ff). However, Anderson's formalization is given in a system which includes the axiom of type-reduction~(\ref{eqn:RM:TRD}), and so is not obviously available once we have rejected this axiom. In his paper \cite{Klement2003aa}, Klement gave a version of the Russell-Myhill paradox which invoked a ``principle of conceivability'' to the effect that ``for every entity, there is at least one sense presenting it as referent'' (\cite{Klement2003aa} p. 305, cf. \S{5} pp. 309~ff). Our Surjectivity Axiom~(\ref{eqn:RM:SM}) is just another expression of  Klement's principle of conceivability. However, in that paper, Klement worked in a system which collapsed the concept-object distinction, so that concepts were a particular species of object.\footnote{See axiom ``PCE'' on \cite{Klement2003aa} p. 305. Another way of formalizing the system of Klement~\cite{Klement2003aa} might be to regard it simply as an untyped system, where there is no distinction between concepts and objects.} Since we want to work within Church's intensional logic, which is a typed system, our formalization has to proceed slightly differently. That said, one can view the formalization given above as the minimal way to modify the formalization of Klement~\cite{Klement2003aa} into the framework of what we're calling the core of Church's system~(\ref{eqn:RM:coresystem}). In particular, while Klement's argument postulated that concepts were a particular species of object, the Senses are Objects Axiom~(\ref{eqn:RM:SO}) postulates that senses are a particular species of object.\footnote{It also bears mentioning that Church, Anderson, and Klement additionally  considered formalizations of the Russell-Myhill paradox within an alternative framework of intensional logic that goes under the heading of ``Russellian intensional logic'' (\cite{Church1984aa}, \cite{Anderson1986aa}, \cite{Klement2002aa} pp. 175~ff). This is the general framework which Klement employed in his widely-read \cite{Klement2005ab}. Since this framework was designed to be an alternative to what we're calling ``Church's intensional logic,'' we have not made use of this in our formalization. By the same token, it is beyond the scope of this paper to say whether anything like a predicative response is available in this alternative framework. To say anything definitive would require at least another lengthy consistency proof like that we offer in \S\ref{sec:RM:05}. If it turned out that nothing like a predicative response was available in this alternative setting, this might well indicate a certain lack of robustness to the predicative response offered in this paper.}

Let's now turn to describing extant responses to this formalized version of the Russell-Myhill paradox. Of course, the specific formalization given above is new to this paper; hence, it is not as if previous authors have explicitly addressed this specific rendition of the paradox. However, anyone who has constructed a model of the core of Church's system~(\ref{eqn:RM:coresystem}) has found some way to avoid this paradox, and so we can ask how these consistent formal systems evade the formalized version of the Russell-Myhill paradox. By Kaplan's construction described in the last section (circa equation~(\ref{eqn:kapalandadsfas})), we can view the standard models of possible worlds semantics as models of Church's core system, and it turns out that the Surjectivity Axiom always comes out true on these models. For, on these models, intensions are just functions from possible worlds to extensions, and so for any extension one can consider the ``constant'' intension that picks out that extension at any world. The models traditionally used in possible worlds semantics are so-called ``standard'' models in which the type~$ab$ is interpreted as the set of all functions from entities of type~$a$ to entities of type~$b$, as judged by the ambient set-theoretic metatheory; hence, the type-theoretic version of Cantor's Theorem~(\ref{eqn:RM:typetheoreticRP}) also comes out true on these models. The tradition of possible worlds semantics then avoids the formalized version of the Russell-Myhill paradox by either rejecting the Propositions as Fine-Grained as Objects Axiom~(\ref{eqn:RM:derepropsaxiom}) or the Senses are Objects Axiom~(\ref{eqn:RM:SO}). For, if there are fewer sets of worlds than there are objects in the worlds, then of course there is no injection from objects to propositions, and so the Propositions as Fine-Grained as Objects Axiom~(\ref{eqn:RM:derepropsaxiom}) comes out false. But if there are more sets of worlds than there are objects in the worlds, then there will be properly more functions from worlds to collections of propositions than there are objects, and hence in this case the Senses are Objects Axiom~(\ref{eqn:RM:SO}) comes out false.\footnote{More formally, we suppose that the types are assigned to sets as follows, wherein~$E$ and~$W$ are fixed sets, corresponding to the objects and the worlds respectively (cf.~(\ref{eqn:RM:usualseq})):
\begin{equation}
D_e = E, \hspace{5mm} D_t=\{0,1\}, \hspace{5mm} D_{ab}=D_b^{D_a} = \{f: D_a\rightarrow D_b\}, \hspace{5mm} D_{a^{\prime}} = D_a^W = \{f: W\rightarrow D_a\}
\end{equation}
Suppose that we are working in a set-theoretic metatheory where as usual $\left|X\right|$ is used to denote the cardinality of the set $X$. Then either~$\left|D_{t^{\prime}}\right|<  \left|D_e\right|$ or not. If so, then there is no injection from~$D_e$ to~$D_{t^{\prime}}$ and the Propositions as Fine-Grained as Objects Axiom~(\ref{eqn:RM:derepropsaxiom}) comes out false. Suppose alternatively that~$\left|D_{t^{\prime}}\right|\geq  \left|D_e\right|$. Since we're working in a set-theoretic metatheory, we can then appeal to Cantor's theorem and basic facts about cardinality to obtain that~$\left| D_{(t^{\prime}t)^{\prime}}\right| \geq \left|D_{t^{\prime}t}\right| >\left|D_{t^{\prime}}\right| \geq \left|D_e\right|$. Then it is not the case that~$D_{(t^{\prime}t)^{\prime}}$ (or anything bijective with it) is a subset of~$D_e$ and hence the Senses are Objects Axiom~(\ref{eqn:RM:SO}) comes out false for the specific type of~$a=(t^{\prime}t)$.}

The work of Anderson (cf. \cite{Anderson1980aa}, \cite{Anderson1984aa} pp. 371 ff) represents a distinct response to the Russell-Myhill paradox, on which one rethinks certain elements of Church's Core System~(\ref{eqn:RM:coresystem}). Anderson's basic idea was to modify Church's system so that there was not a single presentation relation~$\Delta$, but rather a series of presentation relations~$\Delta^{(1)}, \Delta^{(2)}, \ldots$. Let's call~$\Delta^{(n)}$ the~\emph{~$n$-th order presentation function}. On this view, one modifies the axioms of Church's Core System~(\ref{eqn:RM:coresystem}) so that there is one of these axioms for each of the~$n$-th order presentation functions. If one wants to present a formalization of the Russell-Myhill paradox, one needs to specialize it to some specific~$n$-th order presentation function. On this conception, it's natural to think that the analogues of the Surjectivity Axiom would be false. While it might be true in certain models that every higher-order entity was~$n$-th order presented for some~$n$, there might not be a single~$n$ which did this for each higher-order entity. Anderson's response to the Russell-Myhill paradox then parallels the ``typed'' responses to the liar paradox (cf. \cite{Anderson1984aa} p. 376).

This is not the place to argue against the various responses to the formalized version of the Russell-Myhill paradox which we have distilled from the writings of Kaplan and Anderson. For one, it seems likely that such an adjudication would ultimately proceed by reference to  larger considerations like the ability of each of the resulting systems to interpret categorical grammar or to provide a satisfactory semantics for belief attributions. Further, before trying to resolve the paradox in favor of one of these responses, it's important to understand whether we have actually exhausted the entirety of the solution space to the paradox. It seems that the predicativity response has been neglected in the extant literature on the Russell-Myhill paradox, and a chief aim of this paper, which we begin on in earnest in the next section, is to describe the general shape of a plausible predicative response to the Russell-Myhill paradox of propositions.\footnote{Obviously other approaches have been neglected as well: for instance, I know of no extant approaches to the Russell-Myhill paradox which adopt the perspective of non-classical logic.}

\section{The Predicative Response to Russell-Myhill}\label{sec:RM:04}

In the last pages we have surveyed how various constructions of models of Church's core system~(\ref{eqn:RM:coresystem}) avoid the formalized version of the Russell-Myhill paradox (\ref{eqn:RM:formalized}). In the extant literature there seems to have been no attempt to respond to this paradox by rejecting the type-theoretic version of Cantor's Theorem~(\ref{eqn:RM:typetheoreticRP}). Part of the reason for this might be that Cantor's theorem is, well, a theorem. So there might be great pressure to not reject it. But this pressure will only be so great as the strength and plausibility of the axioms from which the theorem is derived. The axioms tacit in the derivation given above of Cantor's theorem are instances of the so-called \emph{comprehension schema}. However, there is a long tradition of \emph{predicative mathematics}, stemming from Poincar\'e and Weyl and represented in our day by the likes of Feferman (\cite{Feferman1964aa,Feferman2005ab}), which proceeds by systematically restricting the comprehension schema. Our aim in what follows is simply to set out and examine a predicative response to the Russell-Myhill paradox of propositions.

Intuitively, the comprehension schema is a mechanism for converting formulas to higher-order entities. In the setting of type-theory, it's most expedient to adopt a version which says that any functional formula determines a higher-order function:
\begin{myenumerate}
\item (Typed Comprehension Schema). The \emph{typed comprehension schema} consists of all the axioms
\vspace{-2mm}\[\forall \; z_1, \ldots, z_k \; [[\forall \; x \; \exists ! \; y \; \varphi(x, y, z_1, \ldots, z_k)]\rightarrow \exists \; h \; [\forall \; x \; \varphi(x, h(x), z_1, \ldots, z_k)]]\]
where~$\varphi(x, y, z_1, \ldots, z_k)$ is a formula with all free variables displayed and with free variable~$x$ of type~$a$,~$y$ of type~$b$, while~$h$ is a variable of type~$ab$ that does not appear free in~$\varphi$.\label{eqn:RM:compschem}\index{Typed Comprehension Schema~(\ref{eqn:RM:compschem})}
\end{myenumerate}
In this schema,  ``$\exists ! x \; \theta(x)$'' is simply the standard abbreviation expressive of uniqueness:  ``$\exists \; x \; (\theta(x) \; \& \; \forall \; z \; (\theta(z)\rightarrow z=x))$''. A special case of this is the following, in which there is no requirement that the formula in question be functional in nature:
\begin{myenumerate}
\item (Concept Comprehension Schema). The \emph{concept comprehension schema} consists of all the axioms~$\forall \; z_1, \ldots, z_k \; \exists \; h \; \forall \; x \;  (h(x)=1 \leftrightarrow \psi(x, z_1, \ldots, z_k))$, where~$\psi(x, z_1, \ldots, z_k)$ is a formula with all free variables displayed and with free variable~$x$ of type~$a$, while~$h$ is a variable of type~$at$ which does not appear free in~$\psi$.\label{eqn:RM:compschem0}\index{Concept Comprehension Schema (\ref{eqn:RM:compschem0})}
\end{myenumerate}
To derive this schema from the Typed Comprehension Schema~(\ref{eqn:RM:compschem}), one defines the formula~$\varphi(x,y,z_1, \ldots, z_k)$ to be
\begin{myequation}\label{eqn:defntrickfromgoing}
 [(\psi(x, z_1, \ldots, z_k) \wedge y=1) \vee  (\neg \psi(x, z_1, \ldots, z_k) \wedge y=0)]
\end{myequation}
The reason for wanting the Typed Comprehension Schema as opposed to the mere Concept Comprehension Schema~(\ref{eqn:RM:compschem0}) is that one wants a way to e.g. go from the entity~$q$ of type~$a$ to the constant function~$f_q(x)=q$ of type~$aa$. Since the Concept Comprehension Schema~(\ref{eqn:RM:compschem0}) only delivers entities~$h$ of type~$at$, it itself cannot do this. Finally, there is a natural principle which generalizes rather than specializes the Typed Comprehension Schema. In particular, if one removes the uniqueness clause from the antecedent of this schema, then this becomes a version of the axiom of choice:
\begin{myenumerate}
\item (Typed Choice Schema). The \emph{typed choice schema} consists of all the axioms
\[\forall \; z_1, \ldots, z_k \; [[\forall \; x \; \exists \; y \; \varphi(x, y, z_1, \ldots, z_k)]\rightarrow \exists \; h \; [\forall \; x \; \varphi(x, h(x), z_1, \ldots, z_k)]]\]
where~$\varphi(x, y, z_1, \ldots, z_k)$ is a formula with all free variables displayed and with free variable~$x$ of type~$a$,~$y$ of type~$b$, while~$h$ is a variable of type~$ab$ which does not appear free in~$\varphi$.\label{eqn:RM:compschem3}\index{Typed Choice Schema~(\ref{eqn:RM:compschem3})}
\end{myenumerate}
The Typed Choice Schema~(\ref{eqn:RM:compschem3}) trivially implies the Typed Comprehension Schema~(\ref{eqn:RM:compschem}), which as we remarked above implies the Concept Comprehension Schema~(\ref{eqn:RM:compschem0}).

The predicative response to the Russell-Myhill paradox restricts the Typed Comprehension Schema (\ref{eqn:RM:compschem}) by imposing constraints on the kinds of higher-order quantifiers and higher-order parameters that can occur in the formula. These restrictions are formulated in terms of the notion of the degree of a type:
\begin{myenumerate}\index{Degree of Type (\ref{eqn:RM:degtypinitial})}
\item (Degree of Type) The \emph{degree} of a type is a positive natural number which is defined recursively as follows: \label{eqn:RM:degtypinitial}\vspace{-2mm}
\[\|e\|=\|t\| =1, \hspace{10mm} \|a^{\prime}\| = \|a\|, \hspace{10mm} \|ab\|=\begin{cases}
\|a\|+1      & \text{if~$\|a\|\geq \|b\|$}, \\
\|b\|      & \text{if~$\|a\|<\|b\|$}.
\end{cases}\]
\end{myenumerate}
To illustrate this last clause, note that $\|e(et)\|=\|t(et)\|=2$ and $\|(et)e\|=\|(et)t\|=3$, while $\|a\|<\|ab\|$ and $\|b\|\leq \|ab\|$ for all types $a,b$. Intuitively, the idea is that degree goes up when the entities of type~$ab$ are genuinely of higher order than those entities of type~$b$. For instance, suppose that~$a=t$ so that that there are only two entities of type~$a$, namely the two truth-values~$0$ and~$1$, and suppose that~$b=et$, so that there are many entities of type~$b$, namely as many as there are concepts. Then entities of type~$ab$ are functions from~$\{0,1\}$ to concepts, and so are really just another way of talking about pairs of concepts. Hence, quantifying over entities of type~$ab$ should involve no more higher-order quantification than quantifying over concepts, and so the degree of~$ab$ should be the same as the degree of~$b$ in this case. 

The predicative comprehension schema may then be defined as follows:
\begin{myenumerate}
\item (Predicative Typed Comprehension Schema). The \emph{predicative typed comprehension schema} consists of all the axioms \vspace{-2mm}
\[\forall \; z_1, \ldots, z_k \; [[\forall \; x \; \exists ! \; y \; \varphi(x, y, z_1, \ldots, z_k)]\rightarrow \exists \; h \; [\forall \; x \; \varphi(x, h(x), z_1, \ldots, z_k)]]\]
where~$\varphi(x, y, z_1, \ldots, z_k)$ is a formula with all free variables displayed and with free variable~$x$ of type~$a$,~$y$ of type~$b$, while~$h$ is a variable of type~$ab$ which does not appear free in~$\varphi$, and in addition variable~$z_i$ has type~$c_i$ with~$\|c_i\|\leq \|ab\|$ and all the bound variables in~$\varphi(x,y,z_1, \ldots, z_k)$ have type~$c$ with~$\|c\|< \|ab\|$.
\label{eqn:RM:predcompschem}\index{Predicative Typed Comprehension Schema (\ref{eqn:RM:predcompschem})}
\end{myenumerate}
Using the same trick as above in equation~(\ref{eqn:defntrickfromgoing}), it's easy to see that the Predicative Typed Comprehension Schema~(\ref{eqn:RM:predcompschem}) implies a version of the Concept Comprehension Schema in which there are the same restrictions on the parameters and bound variables appearing in the formula. For the sake of completeness, we state this version here:
\begin{myenumerate}
\item (Predicative Concept Comprehension Schema). The \emph{predicative concept comprehension schema} consists of all the axioms
\[\forall \; z_1, \ldots, z_k \; \exists \; h \; \forall \; x \;  (h(x)=1 \leftrightarrow \psi(x, z_1, \ldots, z_k))\]
where~$\psi(x, z_1, \ldots, z_k)$ is a formula with all free variables displayed and with free variable~$x$ of type~$a$, while~$h$ is a variable of type~$at$ which does not appear free in~$\psi$, and in addition variable~$z_i$ has type~$c_i$ with~$\|c_i\|\leq \|a\|+1$ and all the bound variables in~$\psi(x,z_1, \ldots, z_k)$ have type~$c$ with~$\|c\|< \|a\|+1$.
\label{eqn:RM:predcompschem0}\index{Predicative Concept Comprehension Schema (\ref{eqn:RM:predcompschem0})}
\end{myenumerate}
Finally, one has the predicative version of the Typed Choice Schema~(\ref{eqn:RM:compschem3}):
\begin{myenumerate}
\item (Predicative Typed Choice Schema). The \emph{predicative typed choice schema} consists of all the axioms \vspace{-1mm}
\[\forall \; z_1, \ldots, z_k \; [[\forall \; x \; \exists \; y \; \varphi(x, y, z_1, \ldots, z_k)]\rightarrow \exists \; h \; [\forall \; x \; \varphi(x, h(x), z_1, \ldots, z_k)]]\]
where~$\varphi(x, y, z_1, \ldots, z_k)$ is a formula with all free variables displayed and with free variable~$x$ of type~$a$,~$y$ of type~$b$, while~$h$ is a variable of type~$ab$ which does not appear free in~$\varphi$, and in addition variable~$z_i$ has type~$c_i$ with~$\|c_i\|\leq \|ab\|$ and all the bound variables in~$\varphi(x,y,z_1, \ldots, z_k)$ have type~$c$ with~$\|c\|< \|ab\|$.\label{eqn:RM:compschem4}\index{Predicative Typed Choice Schema (\ref{eqn:RM:compschem4})}
\end{myenumerate}
Again, it's easy to see that the Predicative Typed Choice Schema~(\ref{eqn:RM:compschem4}) is the deductively strongest schema, so that it implies the Predicative Typed Comprehension Schema~(\ref{eqn:RM:predcompschem}), which in turn implies the Predicative Concept Comprehension Schema~(\ref{eqn:RM:predcompschem0}).

Let's illustrate these predicative schemata by reference to the derivation of the type-theoretic version of the Cantor's Theorem~(\ref{eqn:RM:typetheoreticRP}). In particular, recall the following equation where we defined the diagonal function~$d$ of type~$at$:
\begin{myequation}\label{eqn:digaonaliam}
d(x) = 1 \Longleftrightarrow \exists \; f \; (\iota(f) = x \; \& \; f(x)=0) \tag{\ref{eqn:RM:eqndai}}
\end{myequation}
This is a non-predicative instance of the concept comprehension schema. For, while the defining formula of~$d$ has free variable~$x$ of type~$a$ with degree~$\|a\|$, this formula at the same time contains a bound variable~$f$ of type~$at$. But by consulting the definition of the degree of a type~(\ref{eqn:RM:degtypinitial}), we see that~$\|at\|\geq \|a\|+1$, so that the defining equation~(\ref{eqn:RM:eqndai}) of the diagonal function~$d$ is not an instance of the Predicative Concept Comprehension Schema~(\ref{eqn:RM:predcompschem0}), even though it is an instance of the more general Concept Comprehension Schema~(\ref{eqn:RM:compschem0}). Hence, this  illustrates how if one accepts only the Predicative Typed Comprehension Schema~(\ref{eqn:RM:predcompschem}), the traditional proof of the type-theoretic version of Cantor's Theorem~(\ref{eqn:RM:typetheoreticRP}) is blocked.

A similar elementary observation can be used to illustrate the motivation for the constraints on the parameters in the predicative variants of the comprehension schemas. In particular, if we did not have these restrictions, then we could again show that the diagonal function~$d$ from~(\ref{eqn:digaonaliam}) was in fact a higher-order entity. In this paragraph let us fix a type~$a$ and let us reserve~$\alpha$ and subscripted versions thereof for entities of type~$aa$. Then consider the admittedly uninteresting-appearing formula~$\theta(x,y,z)\equiv y=z$, wherein~$x,y,z$ are of type~$a$. Let~$q$ be an entity of type~$a$. Then trivially we have~$\forall \; x \; \exists \; ! \; y \; \theta(x,y,q)$. Then by the Predicative Typed Comprehension Schema~(\ref{eqn:RM:predcompschem}), there is a function~$\alpha_q$ of type~$aa$ such that~$\alpha_q(x)=q$ for all~$x$ of type~$a$. Let~$\Theta(q,\alpha)\equiv \forall \; x \; \alpha(x)=q$, where again~$x$ has type~$a$ and~$\alpha$ has type~$aa$. Then by the arguments given so far in this paragraph, we have that~$\forall \; q \; \exists \; !  \; \alpha \; \Theta(q,\alpha)$. Hence again by the Predicative Typed Comprehension Schema~(\ref{eqn:RM:predcompschem}), there is a function~$\mathcal{C}$ of type~$a(aa)$ (wherein~$\mathcal{C}$ reminds us of the word ``constant'') such that for all~$q$ of type~$a$, we have~$\mathcal{C}(q)=\alpha_q$. Now, again by the Predicative Concept Comprehension Schema~(\ref{eqn:RM:predcompschem0}), consider the following ``higher-order diagonal''~$\mathcal{D}$ function of type~$(aa)t$:
\begin{myequation}\label{eqn:defn:dqbig}
\mathcal{D}(\alpha)=1 \Longleftrightarrow [\exists \; q \; \alpha = \mathcal{C}(q) \; \& \;  \exists \; f \; (\iota(f)=\alpha(q) \; \& \; f(\alpha(q))=0)]
\end{myequation}
Intuitively,~$\mathcal{D}$ is picking out those constant functions~$\alpha_q$ such that~$d(q)=1$ (where~$d$ is the diagonal function from equation~(\ref{eqn:RM:eqndai})). Now,~$\mathcal{D}$ has type~$(aa)t$ with degree~$\|(aa)t\|>\|a\|+1$. If, contrary to fact, there were no restrictions on parameters in the Predicative Concept Comprehension Schema~(\ref{eqn:RM:predcompschem0}), then we could use~$\mathcal{D}, \mathcal{C}$ to define a higher-order entity~$\widetilde{d}$ of type~$at$ as follows:
\begin{myequation}\label{eqn:defnwided}
\widetilde{d}(\widetilde{q})=1 \Longleftrightarrow \mathcal{D}(\mathcal{C}(\widetilde{q}))=1
\end{myequation}
Then one can verify that
\begin{myequation}\label{eqn:whatwewnattoshceck}
\widetilde{d}(\widetilde{q})=1 \Longleftrightarrow \exists \; f \; (\iota(f) = \widetilde{q} \; \& \; f(\widetilde{q})=0)
\end{myequation}
from which it follows that we have shown that the diagonal function~$d$~(\ref{eqn:RM:eqndai}) again exists as a higher-order entity. This is why it is necessary to include restrictions on parameters in the predicative versions of the comprehension schema.

The predicativity response to the Russell-Myhill paradox has two great burdens. First, it is necessary to show that the Predicative Typed Comprehension Schema~(\ref{eqn:RM:predcompschem}) is consistent with the remaining axioms from the formalized version of the Russell-Myhill paradox~(\ref{eqn:RM:formalized}). In the next section, we discharge this burden by proving:
\begin{myenumerate}
\item (Predicative Consistency Theorem). The following formal system is consistent: the core of Church's system~(\ref{eqn:RM:coresystem}), the Surjectivity Axiom~(\ref{eqn:RM:SM}), the Senses are Objects Axiom~(\ref{eqn:RM:SO}),  the Propositions as Fine-Grained as Objects Axiom~(\ref{eqn:RM:derepropsaxiom}), and the Predicative Typed Choice Schema~(\ref{eqn:RM:compschem4}).\label{eqn:RM:predconthm}\index{Predicative Consistency Theorem (\ref{eqn:RM:predconthm})}
\end{myenumerate}
The second burden of the predicativity response is to say something about what motivates the restriction on the Typed Comprehension Schema~(\ref{eqn:RM:compschem}). Again, intuitively this schema says that every functional formula determines a higher-order entity. If one restricts this, one must say something about \emph{when} and \emph{why} a functional formula determines a higher-order entity. In the remainder of this section, we discuss this more philosophical dimension of the predicativity response to the Russell-Myhill paradox.

Poincar\'e and Weyl were the original predicativists. They drew attention to the fact that higher-order definitions are not in general preserved when one keeps the first-order domain fixed but expands the range of the higher-order quantifiers. So Poincar\'e identifies \emph{predicativity} with definitions that are preserved under such expansions: ``a classification is called \emph{predicative} when it is not changed through the introduction of new elements'' (\cite{Heinzmann1986aa} p. 233, \cite{Poincare1910aa} p. 47, cf. \cite{Chihara1973aa} p. 141). Likewise, in his 1918 \emph{Das Kontinuum}, Weyl draws attention to the fact that if one codes real numbers and continuous functions as certain sets of natural numbers (or rationals), then the definition of the class of continuous functions will contain higher-order quantifiers and thus what counts as a continuous function will depend crucially on the extent and range of the higher-order quantifiers:
\begin{quote}
If we regard the principles of definition as an ``\emph{open}'' system, i.e., if we reserve the right to extend them when necessary by making additions, then in general the question of whether a given function is continuous must also remain \emph{open} [\ldots]' (\cite{Weyl1918} p. 66).
\end{quote}
The dual to Weyl's remark is that if one wants to define notions whose extension does not vary with the extent and range of the higher-order quantifiers, then one should restrict attention to definitions which do not contain higher-order quantifiers.

In my view the best motivation for the predicativity restriction is related to these original thoughts of Poincar\'e and Weyl.\footnote{A distinct set of motivations for predicativity constraints come from the apparent affinity of predicativity with types of constructivism. For more on this complicated aspect of the history of predicativity, see Parsons \cite{Parsons2002}. Another important study of the history of predicativity-like conceptions is Goldfarb's \cite{Goldfarb1988aa} study of Russell. Goldfarb suggests that Russell's reasons for endorsing predicativity-like constraints might be related to having systems in which one can quantify over intensional entities like propositions. By contrast, the motivations given here for predicativity constraints are not intended to have  anything to do with constructivity and are intended to apply with equal force to the quantifiers ranging over intensional entities like propositions as to those ranging over extensional entities like concepts.} In more modern terms, these ideas might be expressed in terms of intuitions about the stability of reference. Suppose that one is using a definite description to refer to an object. If minor variations in empirical conditions cause the object to fail to satisfy this description, then the definite description will be a less than efficacious route to reference. There is a natural generalization of this line of thought to the comprehension schema, where for the sake of simplicity we focus on the Concept Comprehension Schema~(\ref{eqn:RM:compschem0}). The idea is that there is a natural way of seeing each instance of this schema as related to a definite description of a higher-order entity. In particular, consider the following instance of the Concept Comprehension Schema~(\ref{eqn:RM:compschem0}):
\begin{myequation}
\exists \; h \; \forall \; x \; (\varphi(x) \leftrightarrow h(x)=1)
\end{myequation}
One can think about this~$h$ as ``the~$\widetilde{\varphi}$,'' where we define: 
\begin{myequation}\label{eqn:RM:mymy2}
\widetilde{\varphi}(h)\equiv \; [\forall \; x \; (\varphi(x) \leftrightarrow h(x)=1)]
\end{myequation}
If the formula~$\varphi(x)$ contains higher-order quantifiers, then whether a given~$h$ satisfies the description~$\widetilde{\varphi}(h)$ may vary with expansions of the range of the higher-order quantifiers. However, when the formula does not itself contain higher-order quantifiers, then whether something satisfies this description will be stable under expansions of the range of the higher-order quantifiers.

The motivation for the Predicative Typed Comprehension Schema~(\ref{eqn:RM:predcompschem}) can then be seen to derive from the intuition that where one uses a definite description to effect reference to a higher-order entity, this description should be stable under variations of the range of the higher-order quantifiers.  There is nonetheless still a residual philosophical challenge for this predicative response. In particular, it must say something about what the pre-theoretic idea is behind the relevant sense of expansion of the range of the higher-order quantifiers. In my view, the best answer to this is tied to the kinds of positive reasons we can give for the Surjectivity Axiom~(\ref{eqn:RM:SM}). The best positive reason to believe this axiom flows from a conception of what we're trying to model: we're not trying to model higher-order entities as they are in some abstract inaccessible third realm, but we're trying to model higher-order entities insofar as they fall within our referential ken. And it's natural to think that our resources for referring to higher-order entities expands over time just as our resources for referring to concrete objects expands over time.\footnote{Presumably this motivation for the restriction on the quantifiers in the Predicative Typed Comprehension Schema~(\ref{eqn:RM:predcompschem}) likewise motivates the restriction on the higher-order parameters in this schema. In this it's helpful to recall the worked-out example above of the higher-order parameters~$\mathcal{D}$ and~$\mathcal{C}$. As one can see by inspecting equations~(\ref{eqn:defnwided}) and~(\ref{eqn:whatwewnattoshceck}) above, higher-order parameters are able to go proxy for higher-order quantifiers. Given this, if one wants to employ a description featuring a higher-order parameter to stably refer to a lower-order entity, then it's natural to require that this higher-order parameter likewise not shift in extension under variations of the range of the higher-order quantifiers.}

However, it seems safe to say that the same motivations in terms of stability of reference cannot be given for the more general Predicative Typed Choice Schema~(\ref{eqn:RM:compschem4}). Like the Predicative Typed Comprehension Schema~(\ref{eqn:RM:predcompschem}), instances of this schema are conditionals which articulate a sufficient condition for the existence of a higher-order entity. However, unlike the Predicative Typed Comprehension Schema~(\ref{eqn:RM:predcompschem}), it does not seem that the sufficient condition offered by the Predicative Typed Choice Schema~(\ref{eqn:RM:compschem4}) can be conceived of as providing any sort of intension which may serve as a mechanism by which to effect reference to the higher-order entity in question. Thus the Predicative Typed Choice Schema~(\ref{eqn:RM:compschem4}) should not be viewed as following from the predicative viewpoint articulated here, but rather should be viewed as a tool which one can consistently avail oneself of. 

This is important to be clear about because it's only in the presence of some choice-like principle that the Formalized Version of the Russell-Myhill Paradox~(\ref{eqn:RM:formalized}) is actually a deductively valid argument. In particular, in the derivation of the inconsistency in \S\ref{sec:RM:03}, it's easy to see that the move from equation~(\ref{eqn:RM:AC1}) to equation~(\ref{eqn:RM:AC2}) is an instance of some choice-like principle, and one can verify that this move will be covered by the Predicative Typed Choice Schema~(\ref{eqn:RM:compschem4}) (and hence the Typed Choice Schema~(\ref{eqn:RM:compschem3})). If one wanted to be very formal, another axiom that should be officially added to the list of the Formalized Version of the Russell-Myhill Paradox~(\ref{eqn:RM:formalized}) should be the Typed Choice Schema~(\ref{eqn:RM:compschem3}). On this way of putting the matter, the predicative response to the paradox is to deny the type-theoretic version of the Cantor's Theorem~(\ref{eqn:RM:typetheoreticRP}) and to remain ambivalent on the Predicative Typed Choice Schema~(\ref{eqn:RM:compschem4}). Again, by the results of the next section, it's consistent for the predicative response to assume the Predicative Typed Choice Schema~(\ref{eqn:RM:compschem4}), but what we have sought to emphasize in these paragraphs is that the reasons which motivate the predicative restrictions on the comprehension schema don't obviously motivate any instances of either the impredicative or predicative choice schema.

\section{The Consistency of the Predicative Response}\label{sec:RM:05}

In this section, we take up the task of demonstrating the Predicative Consistency Theorem~(\ref{eqn:RM:predconthm}). The reader who is uninterested in this result or merely willing to accept it conditionally might elect to pass on directly to the next section, since the details of this section will not be needed for understanding the subsequent sections of this paper. The most fundamental idea in the proof is to replace the use of the cumulative hierarchy of sets as deployed in the usual model theory of type systems with the constructible hierarchy of sets. Before recalling the definitions of the cumulative hierarchy and the constructible hierarchy of sets, let's then first recall how the usual model theory for type theory proceeds. Let's restrict attention to the \emph{extensional} fragment of the types which contains the type~$e$ for objects, the type~$t$ for truth-values, and which contains the type~$ab$ of functions from entities of type~$a$ to entities of type~$b$ whenever it contains type~$a$ and type~$b$; and moreover, let's consider the \emph{extensional} language which is bereft of the presentation symbols and the intensional application symbols and which contains only the extensional application symbols. 

Models for the extensional fragment of type theory usually begin with an assignment of domains~$D_a$ to each type~$a$. The procedure here is that the type~$e$ of objects is assigned some arbitrary domain~$D_e=E$, the type~$t$ of truth-values is assigned the set~$D_t=\{0,1\}$ of truth-values ($0$ for ``false'' and~$1$ for ``true''), and the type~$ab$ is assigned the domain~$D_{ab}$ of all functions~$f:D_a\rightarrow D_b$, which is sometimes written in exponential notation as~$D_b^{D_a}$ (cf. \cite{Heim1998} p. 28, \cite{Gamut1991aa} pp. 84, 121):
\begin{myequation}\label{eqn:RM:usualseq}
D_e = E, \hspace{10mm} D_t=\{0,1\}, \hspace{10mm} D_{ab}=D_b^{D_a} = \{f: D_a\rightarrow D_b\}
\end{myequation}
To see the connection with the cumulative hierarchy of sets, recall that we can identify sets with their characteristic functions (cf. equation~(\ref{eqn:RM:ide})). Further, recall that 
\begin{myenumerate}
\item \emph{The power set~$P(X)$} of a given set~$X$ is defined to be the set of all the subsets~$Y$ of~$X$, that is,~$P(X)=\{Y: Y\subseteq X\}$.\label{eqn:RM:defnpower}\index{Powerset~$P(X)$ (\ref{defn:powerset}) }\label{defn:powerset}
\end{myenumerate}
Hence, we can identify the sets in~$P(D_a)$ with the functions in~$D_{at}$. For the moment, let's write~$D_{at}\approx P(D_a)$ as a shorthand for this identification. Iterating this, we can build the following sequence in a very straightforward manner:
\begin{myequation}
D_e= E, \hspace{5mm} D_{et}\approx P(E), \hspace{5mm} D_{(et)t} \approx P(P(E)), \hspace{5mm} D_{((et)t)t} \approx P( P(P(E))), \ldots
\end{myequation}
From this perspective, the usual model theory for the extensional theory of types is closely related to iterations of the powerset operator.

This kind of sequence is of course also built into the standard conception of the set-theoretic universe, namely the cumulative hierarchy. In particular, the axioms of set theory guarantee that the universe~$V$ of sets is identical to the union of the following sequence of sets~$V_{\alpha}$, where~$\alpha$ is an ordinal (cf. \cite{Kunen1980} p.~95, \cite{Jech2003} p. 64, \cite{Hrbacek1999aa} p. 257):\index{Cumulative Hierarchy~$V_{\alpha}$ (\ref{eqn:ddfn:cummau})}
\begin{myequation}\label{eqn:ddfn:cummau}
 V_0=\emptyset, \hspace{10mm} V_{\alpha+1} = P(V_{\alpha}), \hspace{10mm} V_{\alpha}=\bigcup_{\beta<\alpha} V_{\beta}, \; \alpha \mbox{ limit}
 \end{myequation}
One of G\"odel's many important innovations in set theory was the definition of the constructible hierarchy of sets. The definition of this hierarchy is identical to the definition of the cumulative hierarchy except at the successor steps, where instead of looking at the full powerset of the previous step, one looks at a certain class of definable subsets of the previous step. In particular, we define, in contrast to the definition of the powerset~(\ref{eqn:RM:defnpower}) above:
 \begin{myenumerate}
 \item \emph{The collection of definable subsets~$\mathrm{Defn}(X)$} of a given set is defined to be the set of all subsets~$Y$ of~$X$ such that there is a first-order formula in the language of set theory~$\varphi(x,z_1, \ldots, z_n)$ with all free variables displayed and parameters~$q_1, \ldots, q_n$ from~$X$ such that~$Y=\{x\in X: (X,\in)\models \varphi(x, q_1, \ldots, q_n)\}$.\label{eqn:defn:RM:defn}\index{Definable subsets~$\mathrm{Defn}(X)$ of~$X$ (\ref{eqn:defn:RM:defn})}
 \end{myenumerate} 
 For instance~$X$ is in~$\mathrm{Defn}(X)$ since~$X=\{x\in X: (X,\in)\models x=x\}$ and the empty set~$\emptyset$ is in~$\mathrm{Defn}(X)$ since~$\emptyset=\{x\in X: (X,\in)\models x\neq x\}$. G\"odel then defined the \emph{constructible hierarchy} as follows (cf. \cite{Kunen1980} p.~166, \cite{Kunen2011aa} p. 134, \cite{Jech2003} p. 174, \cite{Devlin1984aa} p. 58):\index{Constructible Hierarchy~$L_{\alpha}$ (\ref{eqn:mymylimit22})}
\begin{myequation}\label{eqn:mymylimit22}
 L_0=\emptyset, \hspace{10mm} L_{\alpha+1} = \mathrm{Defn}(L_{\alpha}), \hspace{10mm} L_{\alpha}=\bigcup_{\beta<\alpha} L_{\beta}, \alpha \mbox{ limit}
 \end{myequation}
 and he defined the \emph{constructible universe}~$L$ to be the union of the sets from the constructible hierarchy. G\"odel further showed that the constructible universe models all of the axioms of set theory, so that this is not just another collection of sets but an alternative set theoretic universe.

The key idea in our proof of the Predicative Consistency Theorem~(\ref{eqn:RM:predconthm}) is to assign the types to levels of the  constructible hierarchy. To do this, we need to work with a very specific kind of level of the constructible hierarchy. This kind of level was first defined by Kripke (\cite{Kripke1964aa}) and Platek (\cite{Platek1966aa}), who had the idea that  some initial segments of the constructible hierarchy can't ``tell'' that they are tall, and actually think that they can be ``shrunk'', and this idea was later famously employed by Jensen (\cite{Jensen1972aa}) in his proof of the uniformization theorem. Formally one defines this important notion as follows (cf. \cite{Jensen1972aa} pp. 256-257, \cite{Schindler2010aa} Definition 2.1 p. 619, \cite{Devlin1984aa} p. 156, \cite{Sacks1990} p. 157, Barwise \cite{Barwise1975ab} Definition V.6.1 p. 174, \cite{Kripke1964aa} p. 162, \cite{Walsh2014ac}):\index{Projectum,~$n$-th~$\rho_n(\alpha)$ (\ref{defn:theprojectummmm})}
 \begin{myenumerate}
\item the \emph{\mbox{}$n$-th projectum}~$\rho_n(\alpha)$ is the smallest~$\rho\leq \alpha$ such that there is a~$\utilde{\Sigma}_n^{L_{\alpha}}$-definable injection~$\iota: L_{\alpha}\rightarrow \rho$. \label{defn:theprojectummmm}
\end{myenumerate}
In this,~$\utilde{\Sigma}_n$-definability is first-order definability in the sense at issue in the definition of~$\mathrm{Defn}(X)$ above in equation~(\ref{eqn:defn:RM:defn}), but restricted to first-order formulas which begin with a block of alternating quantifiers of length~$n$ starting with an existential quantifier (and allowing parameters). Further, when not clear from context, one writes~$\utilde{\Sigma}_n^{L_{\alpha}}$ to emphasize that the definability is with respect to the structure~$L_{\alpha}$.

Now we can proceed to a description of our models. Our most important definition is the following definition of an \emph{intensional position}~(\ref{eqn:RM:defn:intensional:position}). The motivation for this name comes in the subsequent definition of an \emph{intensional hierarchy}~(\ref{eqn:RM:defn:intensional:hierarchy}) which is given by a collection of \emph{intensional positions}, which intuitively are ``positions'' for the higher-order variables within the hierarchy.\footnote{The language of ``positions'' is apt because, as one can see from inspection of the below definitions, one intensional position can provide the interpretation of the~$n$-th order quantifiers in one intensional hierarchy but the interpretation of the~$m$-th order quantifiers in another.} While the definition of an intensional position is admittedly complicated, the broader significance of each element of this definition will be borne out by the subsequent discussion in this section:
\begin{myenumerate}
\item An \emph{intensional position}~$\mathfrak{p}$ is a given by a sextuple~$\mathfrak{p}=(\alpha, \ell, \iota, \mathcal{O}, \pi, \nu)$ wherein (i) the ordinal~$\alpha>\omega$ is a limit , (ii)~$\ell$ is a positive natural number such that~$L_{\alpha}$ is a model of the~$\Sigma_{\ell}$-collection schema and the~$\Sigma_{{\ell}-1}$ separation schema, (iii) the ordinal~$\alpha$ has non-trivial~$\ell$-th projectum~$\alpha_0=\rho_{\ell}(\alpha)<\alpha$ with~$\iota:L_{\alpha}\rightarrow \alpha_0$ a witnessing~$\utilde{\Sigma}_{\ell}^{L_{\alpha}}$-definable injection, (iv) the set~$\mathcal{O}$ is a~$\utilde{\Sigma}_{\ell}^{L_{\alpha}}$-definable subset of~$\alpha_0$, (v)~the map~$\pi:\mathcal{O}\dashrightarrow L_{\alpha}$ is a~$\utilde{\Sigma}_{\ell}^{L_{\alpha}}$-definable partial surjective function such that~$\pi\circ \iota$ is the identity on~$L_{\alpha}$ and such that~$\mathcal{O}\setminus \pi^{-1}(L_{\alpha})$ is~$\utilde{\Sigma}_{\ell}^{L_{\alpha}}$-definable and (vi) the definability in each of the previous items is with respect to the parameter~$\nu$ from~$L_{\alpha}$.\label{eqn:RM:defn:intensional:position}\index{Intensional Position~(\ref{eqn:RM:defn:intensional:position})}
\end{myenumerate}
In part~(ii) of this definition, the {\emph{}$\Sigma_n$-collection schema} is the axiom schema which says that if~$\varphi(x,y)$ is a~$\Sigma_n$-formula and for all~$x\in a$ there is~$y$ such that~$\varphi(x,y)$, then there is a set~$b$ such that for all~$x\in a$ there is~$y\in b$ with~$\varphi(x,y)$. In essence,~$\Sigma_n$-collection says that when for everything in an antecedently specified set~$a$ there is a witness to a~$\Sigma_n$-condition, then at least one witness for everything in~$a$ may be bounded or collected together in another set~$b$. The \emph{$\Sigma_n$-separation schema} is simply the separation schema from the ambient set theory restricted to the case of~$\Sigma_n$-formulas: it says that if~$\varphi(x)$ is a~$\Sigma_n$ formula and~$a$ is a set then there is another set~$b$ such that~$z\in b$ iff~$z\in a\wedge \varphi(z)$. In essence,~$\Sigma_n$-separation just says that all the~$\Sigma_n$-subsets of antecedently specified~$a$ set exist. Further, it's worth mentioning that the concept of an ordinal~$\alpha$ being~\emph{$\ell$-admissible} from \cite{Walsh2014ac} is equivalent to conditions~(i)-(ii) of the definition of an intensional position, so that intensional hierarchies are just certain collections of~$\ell$-admissibles for increasing values of~$\ell$. This generalizes the notion of Kripke-Platek set theory since  in the case~$\ell=1$, a structure~$L_{\alpha}$ is~$\ell$-admissible just in case it is a model of this set theory (\cite{Kripke1964aa}, \cite{Platek1966aa}, Devlin \cite{Devlin1984aa} p. 48, p. 36).

Having all this in place, we may now define the notion of an \emph{intensional hierarchy}:
\begin{myenumerate}
\item An \emph{intensional hierarchy}~$D=(\mathfrak{p}_1, \mathfrak{p}_2, \ldots)$ is given by a countable sequence~$\mathfrak{p}_n=(\alpha_n, \ell_n, \iota_n, \mathcal{O}_n, \pi_n, \nu_n)$ of intensional positions such that~(i) for all~$n,m\geq 1$ it is the case that~$\rho_{\ell_n}(\alpha_n)= \rho_{\ell_m}(\alpha_m)=\alpha_0$, and (ii)~the associated sequence of ordinals is strictly increasing:~$\alpha_0 <\alpha_1 <\alpha_2<\cdots <\alpha_n<\alpha_{n+1}<\cdots$.\label{eqn:RM:defn:intensional:hierarchy}\index{Intensional Hierarchy (\ref{eqn:RM:defn:intensional:hierarchy})}
\end{myenumerate}
One example of an intensional hierarchy is related to definite descriptions. Let~$\lambda$ be a cardinal in G\"odel's constructible universe~$L$, and let~$\kappa=\lambda^{+}$ be the next biggest cardinal in~$L$, as judged by~$L$. Further, let~$M_n=\mathrm{dcl}^{L_{\kappa}}_{\Sigma_n}(\lambda \cup \{\lambda\})$ be the sets in~$L_{\kappa}$ that have~$\Sigma_n$-definite descriptions over~$L_{\kappa}$ with parameters from~$\lambda \cup \{\lambda\}$. Then it can be shown that that~$M_n=L_{\alpha_n}$ for some~$\alpha$ with~$\lambda<\alpha_n<\alpha_{n+1}<\kappa$ and that~$\rho_n(\alpha_n)=\lambda$. For more details on the construction described in this paragraph, see \cite{Walsh2014ac}, and in particular the existence theorem.\footnote{It's worth spelling out exactly how one defines~$\mathcal{O}_n$ and~$\pi_n$, by more specific reference to the details of the Existence Theorem of \cite{Walsh2014ac} and in particular to the function $\theta_n$ defined therein. The simplest way is to take~$\mathcal{O}_n=\theta_n^{-1}(L_{\alpha_n})$ and to define~$\pi_n=\theta_n\upharpoonright \mathcal{O}_n$. Since~$\theta_n:\mathcal{F}_n\dashrightarrow L_{\alpha_n}$ is~$\utilde{\Sigma}_n^{L_{\alpha_n}}$-definable and~$\mathcal{F}_n$ is~$\utilde{\Sigma}_1^{L_{\alpha_n}}$-definable,~$\mathcal{O}_n$ will be~$\utilde{\Sigma}_n^{L_{\alpha_n}}$-definable, and the \emph{total} surjective map~$\pi_n: \mathcal{O}_n\rightarrow L_{\alpha_n}$ will be similarly definable. Because it is total, trivially~$\mathcal{O}_n\setminus \pi_n^{-1}(L_{\alpha_n})$ is ~$\utilde{\Sigma}_n^{L_{\alpha_n}}$-definable because it is, well, empty. Since~$\pi_n$ is designed to provide the interpretation of~$\Delta_{a}$ for each type~$a$ of degree~$n$ (cf. subsequent discussion circa equation~(\ref{eqn:RM:defnpres})), clearly this interpretation clashes with the intended interpretation of the presentation functions, on which they would be partial. To reinstitute partiality, choose any subset~$\mathcal{P}_n\subseteq \mathcal{F}_n\setminus \theta^{-1}_n(L_{\alpha_n})$ which is~$\utilde{\Sigma}_n^{L_{\alpha_n}}$-definable and then define~$\mathcal{O}^{\prime}_n=\theta_n^{-1}(L_{\alpha_n}) \cup \mathcal{P}_n$ and define ~$\pi_n^{\prime}=\theta_n\upharpoonright \mathcal{O}^{\prime}_n$, making sure to build the parameters defining~$\mathcal{P}_n$ into~$\nu_n$. For instance, one could take~$\mathcal{P}_n$ to be any finite subset of~$\mathcal{F}_n\setminus \theta^{-1}_n(L_{\alpha_n})$.}

Each intensional hierarchy~$D$ naturally gives rise to a model of Church's core system~(\ref{eqn:RM:coresystem}). In particular, we assign types to domains as follows:
\begin{myequation}\label{eqn:defn:RM:typestodomains}\index{Type assignment~$a\mapsto D_a$ (\ref{eqn:defn:RM:typestodomains})}
D_e = \alpha_0,\hspace{5mm} D_t=\{0,1\},\hspace{5mm} D_{ab}={D_b}^{D_a} \cap L_{\alpha_{\|ab\|}},\hspace{5mm} D_{a^{\prime}} = \mathcal{O}_{\|a\|}
\end{myequation}
In this, recall that~$a\mapsto \|a\|$ is the degree of the type~$a$, as defined in~(\ref{eqn:RM:degtypinitial}). So the parallel to the usual semantics for extensional type theory becomes vivid. In particular, whereas these usual semantics employ the cumulative hierarchy to assign domains to types, here our semantics for our intensional type theory uses the constructible hierarchy to assign domains to types. For instance, instead of assigning~$ab$ the set~$D_b^{D_a}=\{f:D_a\rightarrow D_b\}$, we only assign it those elements of this set which are in the constructible hierarchy at an appropriate level. So we're only putting those higher-order entities of type~$ab$ in the range of the higher-order quantifiers when it has entered a level of the constructible hierarchy which is coordinated with the degree of the type~$ab$. 

Having defined intensional hierarchies, our next goal is to say how to interpret the extensional application symbols~$(f,x)\mapsto f(x)$, the presentation symbols~$\Delta_a$, and the intensional application symbols~$(f^{\prime}, a^{\prime})\mapsto f^{\prime}\langle x^{\prime}\rangle$ on intensional hierarchies. In providing these interpretations, we shall be associating each intensional hierarchy~$D$ to an \emph{intensional structure}~$\mathbb{D}$\index{Intensional Structure~$\mathbb{D}$} augmented by these interpretations. Further, as we go along, we shall also show that various axioms are true on these intensional structures. However, prior to doing this, we need to state the following elementary result about how the domains~$D_a$ of an intensional hierarchy relate to the sets~$L_{\alpha_n}$:
\begin{myenumerate}
\item\label{prop:location:domain} (Proposition on the Location of Domains) For all~$n\geq 1$, both of the following hold:
\begin{itemize}
\item[] (I) for all types~$a$ with~$\|a\|<n$, there is a~$\Sigma_1$-formula in parameter~$\mu_n$ such that~$D_a$ is the unique element of~$L_{\alpha_n}$ which satisfies this formula, wherein~$\mu_{n}$ is defined by~$\mu_{n}= \langle \nu_1, \ldots, \nu_{n}, \alpha_0, \alpha_1, \ldots, \alpha_{n-1}\rangle$.
\item[] (II) for all types~$a$ with~$\|a\|= n$, the set~$D_a$ is a~$\utilde{\Sigma}_{\ell_{n}}$-definable subset of~$L_{\alpha_{n}}$ in parameter~$\mu_{n}$.
\end{itemize}
\end{myenumerate}
\noindent For a proof, see Appendix~1 \S\ref{app1}. This result is important because it tells us that the domain~$D_a$ is a subset of~$L_{\alpha_{\|a\|}}$, so that we can locate the domain~$D_a$ amongst the levels of the constructible hierarchy by calculating the degree of the type~$a$. Further, from this we can deduce the following:
\begin{myenumerate}
\item (Proposition on Domain and Codomain of Projectum Witnesses) For all types~$a$, one has that restriction~$\iota_{\|a\|}\upharpoonright D_a$ has domain~$D_a$ and codomain~$D_{a^{\prime}}$, i.e.~$\iota_{\|a\|}\upharpoonright D_a: D_{a}\rightarrow D_{a^{\prime}}$.\label{eqn:RM:domcodomproj}\index{Proposition on Domain and Codomain of Projectum Witnesses (\ref{eqn:RM:domcodomproj})}
\end{myenumerate}
To see this, let~$n=\|a\|$. By the Proposition on the Location of Domains,~$D_a$ is a subset of the domain~$L_{\alpha_n}$ of the injection~$\iota_n: L_{\alpha_n}\rightarrow \alpha_0$. Hence the restriction notation~$\iota_{\|a\|}\upharpoonright D_a$ makes good sense. Suppose now that~$x$ is a member of~$D_a$ and set~$y=\iota_{\|a\|}(x)$. Since~$\pi_n \circ \iota_n$ is the identity function on~$L_{\alpha_n}$, we have that~$y$ is in the domain of~$\pi_{\|a\|}$, which by definition is a subset of~$\mathcal{O}_{\|a\|}=D_{a^{\prime}}$. This, in any case, is the elementary argument which characterizes the Domain and Codomain of the Projectum Witnesses~(\ref{eqn:RM:domcodomproj}).

Given an intensional hierarchy~$D$, there is a natural interpretation of the presentation symbols~$\Delta_a$ such that the Typed Sense Determines Reference Axiom~(\ref{eqn:RM:sdr2}) is true on the induced intensional structure~$\mathbb{D}$. In particular, the presentation functional~$\Delta_a$ is interpreted on an intensional hierarchy~$D$ as the binary relation on~$D_{a^{\prime}}\times D_a$ defined by 
\begin{myequation}\label{eqn:RM:defnpres}
\Delta_a (f^{\prime}, f) \Longleftrightarrow \pi_{\|a\|}(f^{\prime})=f
\end{myequation}
That is,~$\Delta_a$ is interpreted as the graph of~$\pi_{\|a\|}$ restricted to~$D_{a^{\prime}}\times D_a$. This definition makes good sense. For, suppose that~$f^{\prime}$ is from~$D_{a^{\prime}} = \mathcal{O}_{\|a\|}$ and~$f$ is from~$D_a$. Then by the Proposition on Location of Domains~(\ref{prop:location:domain}), we have that~$f\in D_a\subseteq L_{\alpha_{\|a\|}}$ and by definition~$\pi_{\|a\|}:\mathcal{O}_{\|a\|}\dashrightarrow L_{\alpha_{\|a\|}}$. 

We just showed how to expand an intensional hierarchy~$D$ to an intensional structure~$\mathbb{D}$ which has an interpretation of the presentation symbols. Now let us verify that the Typed Sense Determines Reference Axiom~(\ref{eqn:RM:sdr2}) is true on the intensional structure~$\mathbb{D}$. For the ease of readability, we reproduce this axiom here:
\begin{itemize}
\item[(\ref{eqn:RM:sdr2})] \emph{Typed Sense Determines Reference}:~$(\Delta_a(f^{\prime},f) \; \& \; \Delta_a(f^{\prime},g)) \Longrightarrow f=g$
\end{itemize}
Suppose that~$\Delta_a(f^{\prime},f)$ and~$\Delta_a(f^{\prime},g)$. Then by definition in equation~(\ref{eqn:RM:defnpres}), we have that~$\pi_{\|a\|}(f^{\prime})=f$ and~$\pi_{\|a\|}(f^{\prime})=g$. Since~$\pi_{\|a\|}:\mathcal{O}_{\|a\|}\dashrightarrow L_{\alpha_{\|a\|}}$ is a partial function, it then follows that~$f=g$.  Since the Typed Sense Determines Reference Axiom~(\ref{eqn:RM:sdr2}) comes out true on this interpretation of the presentation symbols, we have that the presentation symbol is functional. Just as when working in the object language of Church's core system~(\ref{eqn:RM:coresystem}), instead of writing~$\Delta_a(f^{\prime},f)$, we shall write~$\Delta_a(f^{\prime})=f$. Likewise, we shall write~$\Delta_a(f^{\prime})\hspace{-1mm}\downarrow$ to indicate that there is~$f$ such that ~$\Delta_a(f^{\prime})=f$ (cf. discussion subsequent to~(\ref{eqn:RM:sdr}) in \S\ref{sec:RM:02}).

It remains to indicate the interpretation of the extensional application symbols~$(f,x)\mapsto f(x)$ and the intensional application symbols~$(f^{\prime}, a^{\prime})\mapsto f^{\prime}\langle x^{\prime}\rangle$. The extensional application symbols are comparatively straightforward: these are interpreted as the function from~$D_{ab}\times D_a$ to~$D_b$ given by the notion of extensional application from the metatheory. This makes sense because, per the definition of~$D_{ab}$ in equation~(\ref{eqn:defn:RM:typestodomains}), every element of~$D_{ab}$ is a function~$f:D_{a}\rightarrow D_b$. It's perhaps worth underscoring that for each pair of types~$a,b$, there is a separate extensional application symbol in the signature of Church's core system~(\ref{eqn:RM:coresystem}). We can usually ignore this since their interpretation is uniform. But, in what follows, if we need to explicitly display the types of an extensional application symbol, we shall write~$\mathrm{e\mbox{-}app}_{ab}(f,x)$ instead of~$f(x)$ for the extensional application symbols.\index{Symbol, Extensional Application, Typed~$\mathrm{e\mbox{-}app}_{ab}(f,x)$}

Likewise, we shall sometimes write~$\mathrm{i\mbox{-}app}_{ab}(f^{\prime},x^{\prime})$ instead of~$f^{\prime}\langle x^{\prime}\rangle$ for the intensional application symbols, again to highlight the fact that there is one of these symbols for each pair of types~$a,b$.\index{Symbol, Intensional Application, Typed~$\mathrm{i\mbox{-}app}_{ab}(f^{\prime},x^{\prime})$} We interpret these symbols on an intensional hierarchy as follows:
\begin{myequation}\label{eqn:RM:defn:intensional}
\mathrm{i\mbox{-}app}_{ab}(f^{\prime},x^{\prime}) = f^{\prime}\langle x^{\prime}\rangle=\iota_{\|b\|} ((\Delta_{ab} f^{\prime}) (\Delta_a x^{\prime}))
\end{myequation}
From what we know about the interpretation of the presentation functions and the result on the Domain and Codomain of the Projectum Witnesses~(\ref{eqn:RM:domcodomproj}), we see that the intensional application function is a partial function~$\mathrm{i\mbox{-}app}_{ab}:D_{(ab)^{\prime}}\times D_{a^{\prime}}\dashrightarrow D_{b^{\prime}}$. As with the discussion of the partial presentation functions, technically in the formal system we shall identify the intensional application function with its graph, which is a ternary relation on~$D_{(ab)^{\prime}}\times D_{a^{\prime}}\times D_{b^{\prime}}$. As with presentation symbols, when we write~$\mathrm{i\mbox{-}app}_{ab}(f^{\prime},x^{\prime})$ all by itself, it is assumed that this is defined. 

Now, let's show that the Typed Composition Axiom~(\ref{eqn:RM:cca2}) comes out true on this interpretation. For the ease of readability, we reproduce this axiom here:
\begin{itemize}
\item[(\ref{eqn:RM:cca2})]  \emph{Typed Composition}:~$[\Delta_{ab}(f^{\prime})=f \; \& \; \Delta_a(x^{\prime})=x] \Longrightarrow \Delta_b(f^{\prime}\langle x^{\prime}\rangle) = f(x)$
 \end{itemize}
Suppose that~$\Delta_{ab}(f^{\prime})=f$ and~$\Delta_a(x^{\prime})=x$. Then by its definition in equation~(\ref{eqn:RM:defn:intensional}), we see that~$f^{\prime}\langle x^{\prime}\rangle$ is defined. Then we may evaluate the term~$\Delta_b(f^{\prime}\langle x^{\prime}\rangle)$ as follows:
\begin{myequation}
\Delta_b (\iota_{\|b\|} ((\Delta_{ab} f^{\prime}) (\Delta_a x^{\prime})))=  \Delta_b (\iota_{\|b\|} (f(x))) =(\pi_{\|b\|} \circ \iota_{\|b\|}) (f(x)) = f(x) 
\end{myequation}
where the last equality follows from the fact that~$\pi_{\|b\|}\circ \iota_{\|b\|}$ is the identity function on the set~$L_{\alpha_{\|b\|}}$ (cf. clause~(v) in the definition of an intensional position~(\ref{eqn:RM:defn:intensional:position})). This is why the Typed Composition Axiom~(\ref{eqn:RM:cca2}) comes out true on intensional structures.

Finally, let's note why the Surjectivity Axiom~(\ref{eqn:RM:SM}) and the Senses are Objects Axiom~(\ref{eqn:RM:SO}) are rather trivially true on intensional structures. As for surjectivity, suppose that~$f$ is an element of domain~$D_a$. Again, by the Proposition on Location of Domains~(\ref{prop:location:domain}), we have that~$D_a$ is a subset of~$L_{\alpha_{\|a\|}}$. Since~$\pi_{\|a\|}:\mathcal{O}_{\|a\|}\dashrightarrow L_{\alpha_{\|a\|}}$ is partial surjective, choose~$f^{\prime}$ from~$\mathcal{O}_{\|a\|}$ such that~$\pi_{\|a\|}(f^{\prime})=f$. Then since we have the identity~$D_{a^{\prime}}=\mathcal{O}_{\|a\|}$ (cf. equation~(\ref{eqn:defn:RM:typestodomains})) and since~$\Delta_a$ is interpreted as the graph of~$\pi_{\|a\|}$ restricted to~$D_{a^{\prime}}\times D_a$ (cf. equation~(\ref{eqn:RM:defnpres})), we have that~$\Delta_a(f^{\prime})=f$. This is why the Surjectivity Axiom~(\ref{eqn:RM:SM})  comes out true on intensional structures. As for the Senses are Objects Axiom~(\ref{eqn:RM:SO}), suppose that~$a$ is a type. By definition, we have the identities~$D_e=\alpha_0$ and~$D_{a^{\prime}}=\mathcal{O}_{\|a\|}$ (cf. equation~(\ref{eqn:defn:RM:typestodomains})), and by the definition of an intensional position we have that~$\mathcal{O}_{\|a\|}\subseteq \alpha_0$  (cf. part~(iv) of (\ref{eqn:RM:defn:intensional:position})). Hence, on intensional structures, it is indeed the case that every sense or intension is identical to an object. 

The Propositions as Fine-Grained as Objects Axiom~(\ref{eqn:RM:derepropsaxiom}) requires the following definition:
\begin{myenumerate}
\item An intensional position~$\mathfrak{p}=(\alpha, \ell, \iota, \mathcal{O}, \pi, \nu)$ is \emph{expressive} if there is an injection~$\chi:\alpha_0\rightarrow \pi^{-1}(\{0,1\})$ whose graph is an element of~$L_{\alpha}$. An intensional hierarchy is \emph{expressive} if each position in it is expressive.\index{Intensional Position, Expressive (\ref{intensionalpsotiionex}) }\label{intensionalpsotiionex}
\end{myenumerate}
There are expressive intensional hierarchies (cf.  the existence theorem in \cite{Walsh2014ac}), and any expressive intensional hierarchy models the Propositions as Fine-Grained as Objects Axiom~(\ref{eqn:RM:derepropsaxiom}). In particular, take the injection~$\chi_1: \alpha_0\rightarrow \pi^{-1}_1(\{0,1\})$. Since~$D_e=\alpha_0$ and~$D_t=\{0,1\}$ and~$D_{t^{\prime}}=\mathcal{O}_{1}$ (cf. equation~(\ref{eqn:defn:RM:typestodomains})), it follows that~$\pi^{-1}_1(\{0,1\})\subseteq \mathcal{O}_{1}=D_{t^{\prime}}$, so that~$\chi_1: D_e\rightarrow D_{t^{\prime}}$ is an injection. Moreover, this injection also maps objects to propositions which present a truth-value. For, note that~$\Delta_{t}(\chi_1(x))$ is defined for each~$x$ from~$D_e$ since~$\chi_1(x)\in \pi^{-1}_1(\{0,1\})$. Hence expressive intensional structures model the Propositions as Fine-Grained as Objects Axiom~(\ref{eqn:RM:derepropsaxiom}) in an interesting way since we may inject objects into propositions that actually succeed in presenting truth-values. 

To finish the proof of the Predicative Consistency Theorem~(\ref{eqn:RM:predconthm}), it remains to establish that intensional structures are indeed models of the predicative versions of comprehension:
\begin{myenumerate}
\item (Theorem on Consistency of Predicative Comprehension) For every intensional hierarchy~$D$~(\ref{eqn:RM:defn:intensional:hierarchy}), the associated intensional structure~$\mathbb{D}$ models each instance of the Predicative Typed Choice Schema~(\ref{eqn:RM:compschem4}) and hence each instance of the Predicative Typed Comprehension Schema~(\ref{eqn:RM:predcompschem}). \label{thM:RM:big}\index{Theorem on Consistency of Predicative Comprehension (\ref{thM:RM:big})}
\end{myenumerate}
The proof of this is completed in Appendix 2 \S\ref{app2}, since it is comparatively technical in nature. But with this, the proof of the Theorem on the Consistency of Predicative Comprehension~(\ref{thM:RM:big}) is finished.

In this section we have described models of certain extensions of Church's Core System~(\ref{eqn:RM:coresystem}), and before closing this section it's worth dwelling on one feature of these models related to the Senses are Objects Axiom~(\ref{eqn:RM:SO}). This axiom requires that there be non-trivial identities between different types. While perhaps obvious, it's worth underscoring how this effected. Formally, one simply makes identity an untyped binary relation in the definition of well-formed formulas, as is not uncommon in many-sorted logics (cf. \cite{Manzano1996} p. 229, \cite{Feferman1968} p. 16~footnote~10). On this approach and hence in the models described in this section, identity is simply interpreted as the usual identity relation from the ambient metatheory. One immediate consequence of this approach is that it is only items of syntax such as variables and terms which have a unique type, whereas elements of a domain of a model can have more than one type. This happens more often than one might initially suspect. For instance, the standard semantics for second-order logic is routinely formalized in a many-sorted setting wherein models are given by a pair $(M,P(M))$ wherein the non-empty set $M$ serves as the interpretation of the first-order variables and wherein its powerset $P(M)$ serves as the interpretation of the second-order variables. But if $M$ is a transitive set such as an ordinal, then $M$ is a subset of $P(M)$ and so the two are not at all disjoint. The fact that any finite ordinal is both a first-order object and a second-order object in the standard model of second-order arithmetic $(\omega, P(\omega))$ has never engendered any confusion. Similarly, while the Senses are Objects Axiom~(\ref{eqn:RM:SO}) may be objectionable on purely philosophical grounds, the non-trivial identities between types inherent in it pose no problems for the model theory of such typed systems.

\section{Church's Other Axiom and Gallin's Intensional Logic}\label{sec:RM:06}

Church included another axiom in his own formulation which we have omitted in our original description of his intensional logic in~\S\ref{sec:RM:02}. The strongest version of this axiom is the following, where recall that the ``downarrow'' notation~$\downarrow$ indicates that the partial function is defined on that value (cf. discussion immediately after~(\ref{eqn:RM:sdr}) in \S\ref{sec:RM:02}):
\begin{myenumerate}
\item \emph{Iterative Axiom}:~$\forall \; f^{\prime}, g^{\prime}\; [\Delta_{ab}(f^{\prime})\hspace{-1mm}\downarrow\neq \Delta_{ab}(g^{\prime})\hspace{-1mm}\downarrow] \rightarrow$ 
\item[] \hspace{30mm}~$[\exists \; x^{\prime}, x \; (\Delta_a(x^{\prime})=x \; \& \; \Delta_b(f^{\prime}\langle x^{\prime}\rangle) \neq \Delta_b(g^{\prime}\langle x^{\prime}\rangle))]$ \label{eqn:LSD16preface}\index{Iterative Axiom (\ref{eqn:LSD16preface})}
\end{myenumerate}
The motivation for this axiom is less obvious and Church said less explicitly on this subject. In my view, the best way to conceive of the motivation is as being expressive of a priority of lower-order senses over higher-order senses. The idea is that a canonical way to discern a difference between the \emph{presentations} of higher-order senses is via a difference at the level of the \emph{presentations} of propositions. So one knows that the sense of ``wise'' presents a different concept than the sense of ``courageous'' in part because one knows that, say, the sense of ``Zeno is wise'' presents the true while the sense of ``Zeno has courage'' presents the false.

In Church's papers, this axiom was rather expressed contrapositively as follows (cf. Church's Axiom 16 \cite{Church1951ab} p.~19, \cite{Klement2002aa} pp.~108-109, \cite{Anderson1980aa} pp. 219, 224~ff):
\begin{myequation}
[(\forall \; x,x^{\prime} (\Delta_a(x^{\prime})=x \Rightarrow \Delta_{b}(f^{\prime}\langle x^{\prime}\rangle) = f(x))) \; \& \; \Delta_{ab}(f^{\prime})\hspace{-1mm}\downarrow] \Longrightarrow \Delta_{ab}(f^{\prime})=f\label{eqn:LSD16preface2}
\end{myequation}
In the presence of the Surjectivity Axiom~(\ref{eqn:RM:SM}) and the other axioms of Church's core system~(\ref{eqn:RM:coresystem}), this version follows deductively from the Iterative Axiom~(\ref{eqn:LSD16preface}). To show this, suppose that the antecedent of~(\ref{eqn:LSD16preface2}) holds but the consequent fails, so that~$\Delta_{ab}(f^{\prime})\hspace{-1mm}\downarrow\neq\hspace{-1mm}f$. By the Surjectivity Axiom~(\ref{eqn:RM:SM}), choose~$g^{\prime}$ of type~$(ab)^{\prime}$ such that~$\Delta_{ab}(g^{\prime})=f$. Then~$\Delta_{ab}(f^{\prime})\hspace{-1mm}\downarrow\neq \Delta_{ab}(g^{\prime})\hspace{-1mm}\downarrow$. Then by the Iterative Axiom~(\ref{eqn:LSD16preface}), we have that there is~$x^{\prime}, x$ such that~$\Delta_a(x^{\prime})=x$ and~$\Delta_b(f^{\prime}\langle x^{\prime}\rangle)\neq \Delta_b(g^{\prime}\langle x^{\prime}\rangle)$. By the Typed Composition Axiom~(\ref{eqn:RM:cca2}), we then have that 
\begin{myequation}
\Delta_b(f^{\prime}\langle x^{\prime}\rangle)\neq \Delta_b(g^{\prime}\langle x^{\prime}\rangle) = (\Delta_{ab}(g^{\prime}))(\Delta_a(x^{\prime})) = f(x)
\end{myequation}
which contradicts the hypothesis that the antecedent of~(\ref{eqn:LSD16preface2}) is satisfied. This is the sense in which the Iterative Axiom~(\ref{eqn:LSD16preface}) generalizes Church's own axiom in~(\ref{eqn:LSD16preface2}).

It turns out that the Iterative Axiom~(\ref{eqn:LSD16preface}) (and thus also~(\ref{eqn:LSD16preface2})) are true on the models that we have constructed in the previous section. This is not an accident but rather follows from some additional resources that one has available in these models, resources which allow one to interpret a fragment of Gallin's intensional logic. In particular, let us expand Church's core system~(\ref{eqn:RM:coresystem}) with a new function symbol~$\nabla_a$\index{Symbol, Representation Function~$\nabla_a$} for each type~$a$, called the \emph{representation function}, which takes entities of type~$a$ and returns an entity of type~$a^{\prime}$. Intuitively, the idea is that the representation~$\nabla_a$ function takes an extension~$f$ of type~$a$ and returns an intension~$f^{\prime}$ of type~$a^{\prime}$ which presents~$f$. More formally we have the following axiom:
\begin{myenumerate}
\item \emph{Representation Axiom}: For each entity~$f$ of type~$a$, one has that~$\nabla_a(f)$ is an entity of type~$a^{\prime}$ such that~$\Delta_a(\nabla_a(f))=f$.\label{eqn:RM:PresentationRepresentation}\index{Representation Axiom (\ref{eqn:RM:PresentationRepresentation})}
\end{myenumerate}
Note that it follows from this that the representation function~$\nabla_a$ is an injection from entities of type~$a$ to entities of type~$a^{\prime}$. For, suppose that~$\nabla_a(f)=\nabla_a(g)$. Then by applying the presentation function to each side and by applying the Representation Axiom~(\ref{eqn:RM:PresentationRepresentation}) one has that 
\begin{myequation}
f = \Delta_a(\nabla_a(f)) = \Delta_a(\nabla_a(g))=g
\end{myequation}
Finally, it's perhaps also worth explicitly mentioning that the Representation Axiom~(\ref{eqn:RM:PresentationRepresentation}) formally implies the Surjectivity Axiom~(\ref{eqn:RM:SM}). For the Surjectivity Axiom~(\ref{eqn:RM:SM}) says that each entity~$f$ of type~$a$ is presented by some intension~$f^{\prime}$ of type~$a^{\prime}$, and the  Representation Axiom~(\ref{eqn:RM:PresentationRepresentation}) actually says that one can select the intension~$f^{\prime}$ to be equal to the representation~$\nabla_a(f)$.

The models which we have constructed in the previous section admit a natural interpretation of the representation function on which the Representation Axiom~(\ref{eqn:RM:PresentationRepresentation}) comes out true. In particular, given an intensional hierarchy~$D$~(\ref{eqn:RM:defn:intensional:hierarchy}), we may interpret the representation function~$\nabla_a$ as the injection~$\iota_{\|a\|}$ which comes built into the intensional hierarchy (where again~$\|\cdot\|$ denotes the degree function~(\ref{eqn:RM:degtypinitial}) on types). The Representation Axiom~(\ref{eqn:RM:PresentationRepresentation}) comes out true on intensional structures simply because an intensional hierarchy was built around the idea that the interpretation of the representation function is a (right) inverse to the interpretation of the presentation functions (cf. clause~(v) in the definition of an intensional position~(\ref{eqn:RM:defn:intensional:position}) as well as the Proposition on Domain and Codomain of Projectum Witnesses~(\ref{eqn:RM:domcodomproj})). Further, as we verify in Appendix~2 \S\ref{app2}, the resulting intensional structure continues to model the predicative comprehension schemata~(\ref{eqn:RM:predcompschem})-(\ref{eqn:RM:compschem4}).

The representation functions are relevant to Church's Iterative Axiom due to another axiom which holds true on the models from the last section, namely:
\begin{myenumerate}
\item \emph{Characterization of Intensional Application}:~$f^{\prime}\langle x^{\prime}\rangle = \nabla_b((\Delta_{ab}(f^{\prime})(\Delta_a(x^{\prime})))$
\label{eqn:RM:charintfuncapp}\index{Characterization of Intensional Application (\ref{eqn:RM:charintfuncapp})}
\end{myenumerate}
\vspace{-3mm} One can easily check that this axiom comes out true on intensional structures by glancing at how intensional application was defined in equation~(\ref{eqn:RM:defn:intensional}). This axiom just brings into the object language what was implicit in our constructions in the previous section. Intuitively, what this axiom is saying is that intensional application of a functional intension to an intension simply consists in figuring out what extension is presented by each, performing extensional application on these referents, and then going back to an intension via the representation function. Given this axiom, intensional application can be defined in terms of extensional application and the representation and presentation functions.

Before turning to the connection between the representation function and Church's Iterative Axiom~(\ref{eqn:LSD16preface}), let's note that the axioms governing the representation function allow us to capture a fragment of Gallin's intensional logic~$IL$ (cf. \cite{Gallin1975aa} Chapter 1). Gallin's work can be seen as an attempt to axiomatize Montague grammar (\cite{Dowty1981aa} Chapters 6-8, \cite{Gamut1991aa} Chapter 6). Montague grammar in turn can be viewed as an attempt to develop a logic motivated by possible worlds semantics in which one can distinguish between the intension and the extension of a given expression, while at the same time not having to actually quantify over possible worlds in the object-language (cf. \cite{Gallin1975aa} p. 58, \cite{Dowty1981aa} p. 161). To this end, Montague articulated a type system-- now familiar to us-- in which for any type~$a$ there was a type~$sa$ which in the standard model theory is interpreted as functions from worlds to entities of type~$a$. Montague then postulated that for every well-formed expression~$f$ of type~$a$ there is an intension~\small~$\widehat{f}$ \normalsize of type~$sa$ and for every well-formed expression~$f^{\prime}$ of type~$sa$ there is an extension~\footnotesize$\widecheck{f^{\prime}}$\normalsize \;of type~$a$.

Gallin's later axiomatization can be seen as an attempt to see what is true on all models described by Montague. Some of the crucial axioms that Gallin set out in this intensional logic were the following (cf.  \cite{Gallin1975aa} p. 19):
\begin{myenumerate}
\item \emph{Axiom A2}:~$\forall \; x, y \; [\widehat{x} = \widehat{y} \Longrightarrow \widehat{f(x)}=\widehat{f(y)}]$
\item \emph{Axiom A3}:~$[\forall \; x \; (\widehat{f(x)}=\widehat{g(x)})]\Longrightarrow \widehat{f}=\widehat{g}]$
\item \emph{Axiom AS6}:~$\widecheck{\;\widehat{f}\;\;}=f$
\end{myenumerate}
If we interpret the~$sa$ as~$a^{\prime}$ and we interpret~$\widehat{f}$ by~$\nabla_a(f)$ when~$f$ is of type~$a$ and we interpret~$\widecheck{f^{\prime}}$ as~$\Delta_a(f^{\prime})$ when~$f^{\prime}$ is of type~$a^{\prime}$, then we can easily deduce these three axioms. Hence, the system of Church's intensional logic expanded with the resources of the representation function interprets a fragment of Gallin's intensional logic. One can see this observation as a partial converse to Kaplan's aforementioned possible-worlds model of Church's core system (cf. circa~(\ref{eqn:kapalandadsfas})). The other axioms of Gallin's intensional logic concern modal notions and lambda-terms. The system developed in this paper will have little to say about them because on the one hand it is not a modal system, and on the other hand lambda-terms, as is well known, have the force of effecting the satisfaction of the full Typed Comprehension Schema~(\ref{eqn:RM:compschem}), which we do not have in our predicative setting.

In addition to their interest in providing for this interpretation of a fragment of Gallin's intensional logic, the axioms pertaining to the representation function are of interest because they deductively entail Church's Iterative Axiom.\footnote{To see this, suppose that the antecedent of the Iterative Axiom~(\ref{eqn:LSD16preface}) held, so that~$\Delta_{ab}(f^{\prime})$ was defined and not equal to~$\Delta_{ab}(g^{\prime})$. Then let~$f=\Delta_{ab}(f^{\prime})$ and let~$g=\Delta_{ab}(g^{\prime})$. Since~$f,g$ are functional entities of type~$ab$, it must be the case that they differ on some value (cf. (\ref{extensionalidentityfunctions})), so that there is an object~$x$ of type~$a$ such that~$f(x)\neq g(x)$. Now consider the representation~$x^{\prime}=\nabla_a(x)$ of this entity~$x$. Then by the Representation Axiom~(\ref{eqn:RM:PresentationRepresentation}), we have that~$x^{\prime}$ presents~$x$, or that~$\Delta_a(x^{\prime})=x$. Now, using the Characterization of Intensional  Application~(\ref{eqn:RM:charintfuncapp}), let us quickly compute~$f^{\prime}\langle x^{\prime}\rangle$ and~$g^{\prime}\langle x^{\prime}\rangle$:
\begin{myeqnarray}
f^{\prime}\langle x^{\prime}\rangle=\nabla_b((\Delta_{ab}(f^{\prime})(\Delta_a(x^{\prime}))) = \nabla_b (f(x)) \\
g^{\prime}\langle x^{\prime}\rangle=\nabla_b((\Delta_{ab}(g^{\prime})(\Delta_a(x^{\prime}))) = \nabla_b (g(x))
\end{myeqnarray}
Now, to finish the verification of the Iterative Axiom~(\ref{eqn:LSD16preface}), suppose for the sake of contradiction that~$\Delta_b(f^{\prime}\langle x^{\prime}\rangle)=\Delta_b(g^{\prime}\langle x^{\prime}\rangle)$. Then by the previous calculations, we see that~$\Delta_b(\nabla_b(f(x)))=\Delta_b(\nabla_b(g(x)))$. By the Representation Axiom~(\ref{eqn:RM:PresentationRepresentation}), it then follows that~$f(x)=g(x)$, contrary to hypothesis. Hence, this is why the axioms pertaining to the representation function deductively imply the Iterative Axiom~(\ref{eqn:LSD16preface}).} But in the literature on Church's intensional logic, the ideas behind the Iterative Axiom~(\ref{eqn:LSD16preface}) and the associated principle~(\ref{eqn:LSD16preface2}) have been criticized by Parsons and Klement (\cite{Parsons2001aa} p. 517, \cite{Klement2010aa} pp. 165-166). As Anderson later put it, the general concern is that the Iterative Axiom~(\ref{eqn:LSD16preface}) is ``really quite at odds with the heuristic ideas'' of Church's intensional logic, namely the formalization of fine-grained meanings (\cite{Anderson1998aa} p. 161). One way to see the nature of this concern is to adopt the richer perspective where we have access to the representation function. For, the axioms governing the representation function have the following consequence:
\begin{myenumerate}
\item \emph{Characterization of Intensional Injectivity}: A function~$f$ of type~$ab$ is injective if and only if for any~$f^{\prime}$ of type~$(ab)^{\prime}$ such that~$\Delta_{ab}(f^{\prime})=f$, one has that \label{eqn:RM:CharacterizationofIntensionalInjectivity}~$f^{\prime}\langle x^{\prime}\rangle=f^{\prime}\langle y^{\prime}\rangle$ implies~$\Delta_a(x^{\prime})=\Delta_a(y^{\prime})$. \index{Characterization of Intensional Injectivity (\ref{eqn:RM:CharacterizationofIntensionalInjectivity})}
\end{myenumerate}
The proof of this characterization is comparatively straightforward and so we relegate it to a footnote.\footnote{First suppose that~$f$ is injective and that~$f^{\prime}$ presents~$f$ and that we have the identity~$f^{\prime}\langle x^{\prime}\rangle=f^{\prime}\langle y^{\prime}\rangle$. By the Characterization of Intensional Application~(\ref{eqn:RM:charintfuncapp}), we then have the identity
\begin{myequation}
\nabla_b((\Delta_{ab}(f^{\prime}))(\Delta_a(x^{\prime}))) = \nabla_b((\Delta_{ab}(f^{\prime}))(\Delta_a(y^{\prime})))
\end{myequation}
But since the representation function~$\nabla_b$ is an injection and since~$f^{\prime}$ presents~$f$, this reduces to the identity~$f(\Delta_a(x^{\prime})) = f(\Delta_a(y^{\prime}))$ and since~$f$ is an injection, we have~$\Delta_a(x^{\prime})=\Delta_a(y^{\prime})$, which is what we wanted to show. This completes the verification of the left-to-right direction of~(\ref{eqn:RM:CharacterizationofIntensionalInjectivity}). 

For the right-to-left direction of~(\ref{eqn:RM:CharacterizationofIntensionalInjectivity}), suppose that~$f$ satisfies the right-hand side of (\ref{eqn:RM:CharacterizationofIntensionalInjectivity}). Suppose for the sake of contradiction that~$f$ is not an injection, so that~$f(x)=f(y)$ but~$x\neq y$. Then let~$x^{\prime}=\nabla_a(x)$ and~$y^{\prime}=\nabla_a(y)$ and~$f^{\prime}=\nabla_{ab}(f)$, so that the Representation Axiom~(\ref{eqn:RM:PresentationRepresentation}) implies that~$x^{\prime}$ presents~$x$ and~$y^{\prime}$ presents~$y$ and~$f^{\prime}$ presents~$f$. Then we can expand the identity~$f(x)=f(y)$ to
\begin{myequation}
(\Delta_{ab} f^{\prime})(\Delta_a(x^{\prime})) = f(x) = f(y) = (\Delta_{ab} f^{\prime})(\Delta_a(y^{\prime}))
\end{myequation}
Then by applying the representation function~$\nabla_b$ to each side and appealing to the Characterization of the Intensional Application~(\ref{eqn:RM:charintfuncapp}), we have that~$f^{\prime}\langle x^{\prime}\rangle = f^{\prime}\langle y^{\prime}\rangle$. Then by the hypothesis that~$f^{\prime}$ satisfies the right-hand side of (\ref{eqn:RM:CharacterizationofIntensionalInjectivity}), we have that~$\Delta_a(x^{\prime})=\Delta_a(y^{\prime})$, and since~$x^{\prime}$ presents~$x$ and~$y^{\prime}$ presents~$y$, we have that~$x=y$, which contradicts the reductio assumption that~$x\neq y$. This completes the argument that the Characterization of Intensional Injectivity~(\ref{eqn:RM:CharacterizationofIntensionalInjectivity}) follows from our axioms governing the representation function.} The Parsons-Klement concern can be expressed thusly: this characterization grates against some natural intuitions that one might have about fine-grained meanings.\footnote{See in particular Klement~\cite{Klement2010aa} pp. 165-166. But this is only one aspect of the concern of Parsons and Klement. First, Parsons was most interested in the interaction of Church's other axiom~(\ref{eqn:LSD16preface2}) with senses which do not present any referent (\cite{Parsons2001aa} p. 517). Second, Klement was also concerned with unintuitive consequences of~(\ref{eqn:LSD16preface2}) related to the \emph{intent}ionality of senses (cf. \cite{Klement2010aa} p. 164). A deeper question raised by the work of Parsons and Klement is what analogues there are of the Typed Comprehension Schema~(\ref{eqn:RM:compschem}) for type~$(ab)^{\prime}$. This is relevant because Parsons and Klement's counterexamples to Church's other axiom~(\ref{eqn:LSD16preface2})  are engendered by combining senses or intensions of type~$(ab)^{\prime}$ together in various ways, a procedure which would be most naturally warranted by a version of the comprehension schema for senses or intensions of type~$(ab)^{\prime}$.} Indeed, take any non-injective function that might occur naturally in language, like ``the father of.'' If we take~$a=e$ and~$b=e$ and think of all the objects as consisting of persons, then this is a non-injective function~$f$ of type~$ab$. Let's further assume for the sake of concreteness that intensions of persons are definite descriptions of some kind: \textsc{the~$\Phi$}, \textsc{the~$\Psi$}, etc. Then the Characterization of Intensional Injectivity~(\ref{eqn:RM:CharacterizationofIntensionalInjectivity}) implies that there is some sense~$\textsc{the father of}$ which presents the father-of function~$f$ and which is such that the intension \textsc{the father of }$\langle$\textsc{the~$\Phi$}$\rangle$ is the same qua intension as \textsc{the father of }$\langle$\textsc{the~$\Psi$}$\rangle$, despite the fact that the person who is the~$\Phi$ is not the same as the person who is the~$\Psi$. 

But consequences like this seem highly counterintuitive: one might rather have thought that if the \textsc{the father of~$\langle$the best xylophone player$\rangle$} is the same qua intension as \textsc{the father of~$\langle$the best yazheng player$\rangle$}, then the best xylophone player is the best yazheng player. (The ``\emph{x}ylophone'' and ``\emph{y}azheng'' are two musical instruments which start with the same letters that are used as the variables in (\ref{eqn:RM:CharacterizationofIntensionalInjectivity})). However, what we now see is that if one accepts the Characterization of Intensional Application~(\ref{eqn:RM:charintfuncapp}), then one must to accept consequences like this. For, on this characterization of intensional application, the intension associated to the \textsc{the father of~$\langle$the best xylophone player$\rangle$} is \emph{not} the definite description which we normally associate to the linguistic expression ``the father of the best xylophone player.'' Rather, on this characterization of intensional application, the intension \textsc{the father of~$\langle$the best xylophone player$\rangle$} is the result of {\it intensionally} applying the intensional functional~\textsc{the father of} to the input of the intension \textsc{the best xylophone player}.

Indeed, on the conception following from the Characterization of Intensional Application~(\ref{eqn:RM:charintfuncapp}), this is done by first by figuring out who the father of the best xylophone player actually is-- perhaps its Ted-- and going and figuring out what the representation of Ted is -- perhaps it is \textsc{the mayor of Montreal}. Now Ted might have two children, Alice and Bob, and it might turn out that Alice is best xylophone player while Bob is the best yazheng player. On this conception, the intension associated to the \textsc{the father of~$\langle$the best xylophone player$\rangle$} is identical to the intension associated to \textsc{the father of~$\langle$the best yazheng player$\rangle$} since both are identical to the intension \textsc{the mayor of Montreal}. But in spite of this identity, the person presented by the intension \textsc{the best xylophone player} is Alice, who is distinct from her sibling Bob, who is presented by the intension \textsc{the best yazheng player}.

Thus the Characterization of Intensional Application~(\ref{eqn:RM:charintfuncapp}) requires us to depart from some of the original ambitions of a fine-grained theory of intensions, on which intensional injectivity would presumably be the rule rather than the exception.\footnote{For instance, some of the systems of Church and Anderson explicitly included axioms for the injectivity of senses of functional expressions. See the axiom designated ``64'' in Church \cite{Church1974aa} p. 151 and Anderson \cite{Anderson1980aa} p. 222.} It is not presently obvious to us whether there is a proof of the Predicative Consistency Theorem~(\ref{eqn:RM:predconthm}) which would produce models which do not validate either the Characterization of Intensional Injectivity~(\ref{eqn:RM:CharacterizationofIntensionalInjectivity}) or the Iterative Axiom~(\ref{eqn:LSD16preface}). If one rejects Church's Axiom of Type Reduction~(\ref{eqn:RM:TRD}) then the most difficult part of any construction of a model of these systems is to provide an interpretation of the intensional application function~$(f^{\prime}, x^{\prime})\mapsto f^{\prime}\langle x^{\prime}\rangle$. All the constructions which we have come up with so far have involved the aforementioned characterization of intensional application~(\ref{eqn:RM:charintfuncapp}) and hence the Iterative Axiom~(\ref{eqn:LSD16preface}). 

\section{Wehmeier and the Problem of Many Non-Extensions}\label{sec:RM:07}

In this paper, we've been primarily concerned with describing the predicative response to the Russell-Myhill paradox. However, as we've seen, predicativity constraints block the normal proof of the type-theoretic version of Cantor's Theorem~(\ref{eqn:RM:typetheoreticRP}). Since this theorem is closely related to Russell's paradox, it's natural to think that there is a connection between the predicative response to the Russell-Myhill paradox of propositions and consistent fragments of the so-called naive conception of set.

To see this connection, let's note the precise way in which we can use axioms introduced thus far to produce violations of the type-theoretic version of Cantor's Theorem~(\ref{eqn:RM:typetheoreticRP}). In particular, let's note why these axioms give us reason to endorse the following principle:
\begin{myenumerate}
\item There is an injection~$\partial$ from entities of type~$et$ to entities of type~$e$ such that for all~$f$ of type~$et$ there is~$f^{\prime}$ of type~$(et)^{\prime}$ satisfying~$\Delta_{et}(f^{\prime})=f$ and~$f^{\prime}=\partial(f)$. \label{eqn:32143214123}
\end{myenumerate}
One proof of this proceeds via the representation function~$\nabla_{et}$ introduced in the previous section (cf. circa~equation~(\ref{eqn:RM:PresentationRepresentation})). For, one can use the Predicative Typed Comprehension Schema~(\ref{eqn:RM:predcompschem}) and the Senses are Objects Axiom~(\ref{eqn:RM:SO}) to define the map~$\partial$ as follows, where~$f$ is a variable of type~$et$,~$x$ is a variable of type~$e$, and~$f^{\prime}$ is a variable of type~$(et)^{\prime}$:
\begin{myequation}\label{eqn:almostdone}
\partial(f) = x \Longleftrightarrow \exists \; f^{\prime} \; (f^{\prime}=\nabla_{et}(f) \; \& \; f^{\prime}=x)
\end{myequation}
While the representation function~$\nabla_{et}$ is a function from entities of type~$et$ to entities of type~$(et)^{\prime}$, the~$\partial$ function is a function from entities of type~$et$ to entities of type~$e$. Hence~$\partial$ is an injection since as we noted in the last section the representation function~$\nabla_{et}$ is an injection.

A second proof of~(\ref{eqn:32143214123}) proceeds by recourse to the Predicative Typed Choice Schema~(\ref{eqn:RM:compschem4}). For, by the Surjectivity Axiom~(\ref{eqn:RM:SM}) and the Senses are Objects Axiom~(\ref{eqn:RM:SO}), one has the following, where again~$f$ is a variable of type~$et$,~$x$ is a variable of type~$e$, and~$f^{\prime}$ is a variable of type~$(et)^{\prime}$:
\begin{myequation}
\forall \; f \; \exists \; x \; [\exists \; f^{\prime} \; \Delta_{et}(f^{\prime})=f \; \& \; f^{\prime}=x]
\end{myequation}
Then by the  Predicative Typed Choice Schema~(\ref{eqn:RM:compschem4}), it follows that there is a function~$\partial$ of type~$(et)e$ such that
\begin{myequation}
\forall \; f \; [\exists \; f^{\prime} \; \Delta_{et}(f^{\prime})=f \; \& \; f^{\prime}=\partial(f)]
\end{myequation}
Then we may argue that~$\partial$ is an injection. For suppose that~$\partial(f)=\partial(g)$. Then by the previous equation, there are~$f^{\prime}, g^{\prime}$ such that~$f^{\prime}=\partial(f)=\partial(g)=g^{\prime}$ and~$\Delta_{et}(f^{\prime})=f$ and~$\Delta_{et}(g^{\prime})=g$. Then since~$f^{\prime}=g^{\prime}$, we have that~$f=\Delta_{et}(f^{\prime})=\Delta_{et}(g^{\prime})=g$, so that the injectivity of~$\partial$ is thereby established.

One of the most traditional versions of the naive conception of set is that found in Frege's \emph{Grundgesetze} (\cite{Frege1893}, \cite{Frege2013aa}). One of the crucial axioms of this system is Basic Law~V, which postulates the existence of a injection from concepts to objects, which we may call the \emph{extension operator}. Now, concepts can be identified with functions of type~$et$ (as we have had numerous occasions to observe in this paper, e.g. circa equation~(\ref{eqn:RM:ide})). What the previous paragraphs then show is that the expansions of Church's intensional logic which we have studied in this paper afford the resources to satisfy one key postulate of Frege's \emph{Grundgesetze}-- namely Basic Law~V-- along with fragments of the comprehension schema, such as the Predicative Typed Comprehension Schema~(\ref{eqn:RM:predcompschem}). This consistency result in and of itself is not new: versions of it were established by Parsons~\cite{Parsons1987a}, Heck~\cite{Heck1996}, and Ferreira-Wehmeier~\cite{Ferreira2002aa}, and it was the focus of parts of our earlier papers \cite{Walsh2012aa}, \cite{Walsh2014ac}. However, the argument of the previous paragraphs is new in that it establishes the existence of a particular {\it species} of extension operator, which we might dub a {\it sense-selecting extension operator} and define formally as follows:
\begin{myenumerate}
\item A \emph{sense-selecting extension operator} is an injection~$\partial$ from entities of type~$et$ to entities of type~$e$ such that for all~$f$ of type~$et$ there is~$f^{\prime}$ of type~$(et)^{\prime}$ satisfying~$\Delta_{et}(f^{\prime})=f$ and~$f^{\prime}=\partial(f)$.\label{eqn:defn:sensebearing}\index{Sense-Selecting Extension Operator (\ref{eqn:defn:sensebearing})}\index{Symbol, Sense-Selecting Extension Operator~$\partial$}
\end{myenumerate}
Part of what is added by looking at Frege's naive conception of set as~embedded within a certain expansion of Church's intensional logic is that we have access to a particular kind of extension operator, one on which the extension of a concept is a sense of that concept.\footnote{However, it should be emphasized that the predicative response, as we have described it above in \S\ref{sec:RM:04}, is not necessarily committed to the existence of a sense-selecting extension operator. For, the two proofs from the above paragraphs used the representation operator~$\nabla_a$ from \S\ref{sec:RM:06} and the Predicative Typed Choice Schema~(\ref{eqn:RM:compschem4}). As stressed in \S\ref{sec:RM:04}, the philosophical motivations for the Predicative Typed Comprehension Schema~(\ref{eqn:RM:predcompschem0}) don't necessarily extend to the Predicative Typed Choice Schema~(\ref{eqn:RM:compschem4}); and it goes without saying that while the representation operator helps bring more of the model construction into the object-language, it too is not necessarily built into the predicative response to the Russell-Myhill paradox. Indeed, it is not even clear to me whether one can derive the existence of sense-selecting extension operators merely from the core of Church's system~(\ref{eqn:RM:coresystem}), the Surjectivity Axiom~(\ref{eqn:RM:SM}), the Senses are Objects Axiom~(\ref{eqn:RM:SO}), the Propositions as Fine-Grained as Objects Axiom~(\ref{eqn:RM:derepropsaxiom}), and the Predicative Typed Comprehension Schema~(\ref{eqn:RM:predcompschem0}). Thus the results of this section are only available to certain natural expansions of the predicative perspective by choice principles or by a representation operator.
} Of course this general kind of maneuver is familiar from the literature on the philosophy of set theory. For instance, the stage axioms of Shoenfield~\cite{Shoenfield1961aa}, \cite{Shoenfield1967aa}, \cite{Shoenfield1977aa} and Boolos~\cite{Boolos1971} constitute an embedding of a fragment of Zermelo-Fraenkel set theory within a theory of ``collections formed in stages.'' This gave Shoenfield and Boolos additional resources by which to respond to the Quinean charge that this set theory was just ``wisdom after paradox,'' or just one of many ad-hoc responses to the paradoxes (\cite{Quine1986aa} p.~403, cf. \cite{Quine1984aa} p.~789, \cite{Quine1960aa} pp.~353-354, \cite{Quine1969aa} p.~5, \cite{Martin1970aa} pp.~111-112).

Similarly, viewing Frege's set theory in the light of Church's intensional logic allows us to respond to a serious objection, due to Wehmeier, with these consistent fragments of the \emph{Grundgesetze}. One way to see one's way towards this objection is to observe that working within the framework of Church's intensional logic, we can show that sense-selecting extension operators have ranges which are indefinitely extensible in the sense of Russell and Dummett. To this end, let us first define some preliminary subset notation:
\begin{myenumerate}
\item \emph{Subset Notation}: If~$\Phi(x)$ is a formula in one free variable~$x$ of type~$a$ and~$f$ is an entity of type~$at$, let's say that~$h\subseteq \Phi$ if~$\forall \; x \; (h(x)=1 \rightarrow \Phi(x))$. Likewise, if~$g$ is also of type~$at$, let's say that~$h\subseteq g$ iff~$\forall \;x \; (h(x)=1 \rightarrow g(x)=1)$, and let us define~$h\subsetneq g$ as~$h\subseteq g \wedge \neg (g\subseteq h)$.\index{Symbol, Subset~$h\subseteq \Phi(x)$ (\ref{eqn:subsetneo})}\label{eqn:subsetneo}
\end{myenumerate}
Then we may define a formal version of indefinite extensibility as follows:
\begin{myenumerate}
\item A formula~$\Phi(x)$ in one free variable~$x$ of type~$a$ is \emph{formally indefinitely extensible} if for each~$h$ of type~$at$ with~$h\subseteq \Phi$ there is~$\widetilde{h}$ of type~$at$ such that~$h\subsetneq \widetilde{h}\subseteq \Phi$.\label{eqn:FormIndefExt}\index{Formally Indefinitely Extensible (\ref{eqn:FormIndefExt})}
\end{myenumerate}
Dummett, following Russell, expressed the idea of indefinite extensibility as follows: ``[a]n indefinitely extensible concept is one such that, if we can form a definite conception of a totality all of whose members fall under that concept, we can, by reference to that totality, characterize a larger totality of all whose members fall under it'' (\cite{Dummett1994aa} p. 22, \cite{Dummett1963aa} pp. 149-150, \cite{Dummett1978aa} pp. 195-196, \cite{Dummett1981aa} p. 533, \cite{Dummett1991aa} p. 316, cf. \cite{Russell1907aa} p. 36, \cite{Russell1973aa} p. 144). If one reads Dummett's use of ``concept'' as any formula~$\Phi(x)$ with a free object variable~$x$ and if one reads his ``definite concept'' as an entity of type~$et$, then there is a comparatively tight match between the formalization in~(\ref{eqn:FormIndefExt}) and Dummett's own formulation of indefinite extensibility.\footnote{That said, there are some differences. First, this formalization provides no insight into how,  if at all,~$\widetilde{h}$ is provided ``by reference'' to~$h$. Second, on our explication of ``definite'', it will follow that the definite concepts are closed under boolean operations such as intersection, union, and complement. If one has the intuition that ``definite concepts'' should be small in some sense, one will resist the claim that definite concepts are closed under complementation. Finally, it should be noted that this general variety of formalization of indefinite extensibility is of course not new: see for instance Shapiro-Wright \cite{Shapiro2006ac} p. 266 and Priest \cite{Priest2013aa} pp. 1264-1265.}

Let's now show that if~$\partial$ is a sense-selecting extension operator~(\ref{eqn:defn:sensebearing}) then the range~$\mathrm{rng}(\partial)$ of this extension operator is formally indefinitely extensible~(\ref{eqn:FormIndefExt}), where of course the range~$\mathrm{rng}(\partial)$ is the following formula with~$x$ a variable of type~$e$ and~$f$ a variable of type~$et$:
\begin{myequation}
(\mathrm{rng}(\partial))(x)\equiv \exists \; f \; \partial(f)=x
\end{myequation}
Fix~$h$ of type~$et$ such that~$h\subseteq \mathrm{rng}(\partial)$. Then one may show the following, which intuitively says that~$\partial$ admits a partial inverse:
\begin{myenumerate}
\item There is a~$\gamma_h$ of type~$e(et)$ such that for every~$g$ of type~$et$ with~$h(\partial(g))=1$, it is the case that~$\gamma_h(\partial(g))=g$.\label{eqn:surjectiveproperty}
\end{myenumerate}
Since the argument for~(\ref{eqn:surjectiveproperty}) is routine, we relegate it to a footnote.\footnote{The first argument for~(\ref{eqn:surjectiveproperty}) employs the representation function. So suppose that the sense-selecting extension operator satisfies $\partial(f)=\nabla_{et}f$ as in equation~(\ref{eqn:almostdone}). Fix a parameter $q$ of type $et$. Then one has the following, wherein $x$ has type~$e$ and $f$ has type~$et$:
\begin{myequation}
\forall \; x \; \exists \; ! \; f \; [(h(x)=0 \; \& \; f=q) \vee (h(x)=1 \; \& \; \nabla_{et}(f)=x)]
\end{myequation}
Then by the Predicative Typed Comprehension Schema~(\ref{eqn:RM:predcompschem}), there is $\gamma_h$ of type $e(et)$ such that
\begin{myequation}
\forall \; x \; [(h(x)=0 \; \& \; \gamma_h(x)=q) \vee (h(x)=1 \; \& \; \nabla_{et}(\gamma_h(x))=x)]
\end{myequation}
To verify equation~(\ref{eqn:surjectiveproperty}), suppose that $h(\partial(g))=1$. Letting $x=\partial(g)$ we have that $h(x)=1$. Then $\nabla_{et}(\gamma_h(x))=x=\partial(g)=\nabla_{et}g$. Then $\gamma_h(x)=g$ and so $\gamma_h(\partial(g))=g$, which is what we wanted to show. The second argument for~(\ref{eqn:surjectiveproperty}) employs the Predicative Typed Choice Schema~(\ref{eqn:RM:compschem4}). Since~$h\subseteq \mathrm{rng}(\partial)$, we have that~$\forall \; x \; [h(x)=1\rightarrow (\exists \; f \; \partial(f)=x)]$. Then the definition of sense-selecting~(\ref{eqn:defn:sensebearing}) implies that~$\forall \; x \; [h(x)=1\rightarrow (\exists \;  f \; \exists \; f^{\prime} \; \Delta_{et}(f^{\prime}) = f \wedge f^{\prime}=x)]$. Trivially we then have~$\forall \; x \; \exists \; f^{\prime} \; [h(x)=1\rightarrow f^{\prime}=x]$. Then we may apply the Predicative Typed Choice Schema~(\ref{eqn:RM:compschem4}) since the parameter~$h$ has type with degree~$2$ and the type~$e(et)^{\prime}$ has degree~$2$. Doing this we get an entity~$\beta_h$ of type~$e(et)^{\prime}$ such that~$\forall \; x \; [h(x)=1\rightarrow \beta_h(x)=x]$. Further, we claim that
\begin{myequation}\label{eqn:mymy34}
\forall \; x \; \exists \; f \; [h(x)=1\rightarrow \Delta_{et}(\beta_h(x))=f]
\end{myequation}
For, if~$h(x)=1$ then~$\partial(f)=x$ for some~$f$ of type~$et$ and hence~$\Delta_{et}(g^{\prime})=f$ and~$g^{\prime}=x$ for some~$g^{\prime}$ of type~$(et)^{\prime}$ by the definition of sense-selecting~(\ref{eqn:defn:sensebearing}). Then~$\beta_h(x)=x=g^{\prime}$ and so~$\Delta_{et}(\beta_h(x))=\Delta_{et}(g^{\prime})=f$. So indeed equation~(\ref{eqn:mymy34}) holds. Further, we may  apply the Predicative Typed Choice Schema~(\ref{eqn:RM:compschem}) to this equation since the parameters~$h,\beta_h$ have types with degree~$2$ and since~$e(et)$ likewise has degree~$2$. Then we obtain~$\gamma_h$ of type~$e(et)$ such that
\begin{myequation}\label{eqn:mymy35}
\forall \; x \; [h(x)=1\rightarrow \Delta_{et}(\beta_h(x))=\gamma_h(x)]
\end{myequation}
Let's now verify equation~(\ref{eqn:surjectiveproperty}). Suppose that~$g$ is of type~$et$ such that~$h(\partial(g))=1$. Let~$x=\partial(g)$, so that~$h(x)=1$. By~$\partial(g)=x$, we obtain~$\Delta_{et}(g^{\prime})=g$ and~$g^{\prime}=x$ for some~$g^{\prime}$ of type~$(et)^{\prime}$ by the definition of sense-selecting~(\ref{eqn:defn:sensebearing}). Further by equation~(\ref{eqn:mymy35}), we have ~$\Delta_{et}(\beta_h(x))=\gamma_h(x)$. Then~$\beta_h(x)=x=g^{\prime}$ and so~$g=\Delta_{et}(g^{\prime})=\Delta_{et}(\beta_h(x))=\gamma_h(x)=\gamma_h(\partial(g))$, which is what we wanted to show.} By the definition of degree in equation~(\ref{eqn:RM:degtypinitial}), note that~$\gamma_h$ has type with degree~$2$. Hence the following formula~$\varphi(x,h,\gamma_h)$, where~$x$ is a variable of type~$e$, contains only parameters of degree~$2$:
\begin{myequation}
\varphi(x,h,\gamma_h) \equiv (h(x)=1 \; \& \; (\gamma_h(x))(x)=0)
\end{myequation}
Then by the Predicative Concept Comprehension Schema~(\ref{eqn:RM:predcompschem0}), there is~$g_h$ of type~$et$ such that
\begin{myequation}
g_h(x) =1 \Longleftrightarrow  (h(x)=1 \; \& \; (\gamma_h(x))(x)=0)
\end{myequation}
Then we claim that~$h(\partial(g_h))=0$. For, suppose not. Then let~$y=\partial(g_h)$ so that~$h(y)=1$. Then by the earlier result~(\ref{eqn:surjectiveproperty}) we have that~$\gamma_h(y)=\gamma_h(\partial(g_h))=g_h$. Then the above equation implies that 
\begin{myequation}
g_h(y) =1 \Longleftrightarrow  (\gamma_h(y))(y)=0 \Longleftrightarrow g_h(y)=0
\end{myequation}
which is a contradiction. Hence indeed we have~$h(\partial(g_h))=0$. Now, let~$p=\partial(g_h)$, so that~$p$ is a parameter of type~$e$ with degree~$1$. Then consider the following formula which has only parameters of degree~$\leq{2}$:
\begin{myequation}
\psi(x,h,p) \equiv (h(x)=1 \vee x=p)
\end{myequation}
Then by the Predicative Concept Comprehension Schema~(\ref{eqn:RM:predcompschem0}), there is~$\widetilde{h}$ of type~$et$ such that
\begin{myequation}
\widetilde{h}(x) =1 \Longleftrightarrow (h(x)=1 \vee x=\partial(g_h))
\end{myequation}
so that in terms of our subset notation (\ref{eqn:subsetneo}), we have~$h\subsetneq \widetilde{h} \subseteq \mathrm{rng}(\partial)$ which completes the verification that the range~$\mathrm{rng}(\partial)$ of a sense-selecting extension operator is formally indefinitely extensible~(\ref{eqn:FormIndefExt}).

Let's call the objects falling within the range~$\mathrm{rng}(\partial)$ of an extension operator~$\partial$ the \emph{extensions}. In this terminology, Wehmeier's observation about his consistent fragments of the \emph{Grundgesetze} was that they required that there were infinitely many non-extensions (cf. \cite{Wehmeier1999} \S{4.2} pp.~326~ff, \cite{Wehmeier2004aa} \S{3} pp.~255~ff). Given the above discussion, one can see now that this follows deductively from the formally indefinite extensibility of the range~$\mathrm{rng}(\partial)$ of an extension operator~$\partial$. For, suppose that there were only finitely many objects which were non-extensions, enumerated as~$q_1, \ldots, q_k$. Then consider the formula~$\theta(x,q_1, \ldots, q_k)$ with free variable~$x$ of type~$e$ and parameters~$q_1, \ldots, q_k$ of type~$e$ and hence degree~$\leq{1}$:
\begin{myequation}
\theta(x,q_1, \ldots, q_k) \equiv (x\neq q_1 \wedge \cdots \wedge x\neq q_k)
\end{myequation}
Then by the Predicative Concept Comprehension Schema~(\ref{eqn:RM:predcompschem0}), there is~$h$ of type~$et$ such that 
\begin{myequation}
h(x)=1 \Longleftrightarrow (x\neq q_1 \wedge \cdots \wedge x\neq q_k)
\end{myequation}
But by hypothesis, we have that~$h$ is coextensive with~$\mathrm{rng}(\partial)$, in that~$h(x)=1$ iff~$x$ is in the range of the extension operator~$\partial$. But then by the formal indefinite extensibility~(\ref{eqn:FormIndefExt}) of~$\mathrm{rng}(\partial)$, there is~$\widetilde{h}$ of type~$et$ such that~$h\subsetneq \widetilde{h} \subseteq \mathrm{rng}(\partial)$, which contradicts that~$h=\mathrm{rng}(\partial)$. Thus the supposition that there were only finitely many non-extensions must have been wrong. Hence, the formal indefinite extensibility of the range of the extension operator requires that there be infinitely many non-extensions.

As Wehmeier notes (cf. \cite{Wehmeier1999} \S{4.2} pp.~326~ff, \cite{Wehmeier2004aa} \S{3} pp.~255~ff),  these considerations suggest an apparent tension between these subsystems of the \emph{Grundgesetze} and at least some renditions of Frege's logicism. For, sometimes logicism is described as the contention that mathematical reasoning is discoverable in every domain of inquiry (cf. (\cite{Demopoulos2005ab} p.~138, cf. \cite{Demopoulos1998aa} \S{VII} p.~496, \cite{Demopoulos1994aa} p.~229). Presumably there are domains of inquiry (such a chemistry and biology) in which there are comparatively few non-extensions (say, finitely many atoms or organisms). In such domains of inquiry there simply isn't ``space enough'' for an extension operator as axiomatized by the predicative fragments of the \emph{Grundgesetze} of the kind considered here. Besides this apparent tension, there is a more general reason to be concerned about the problem of many non-extensions. For, this problem tells us that the presence of an extension operator has consequences for the non-extensions. It's natural to seek an explanation for this-- that is, one seeks an answer to the question: what is it about the extension operator that results in it having consequences for the nature of the non-extensions?

But in the case where the extension operator is a sense-selecting extension operator, it seems that there is a natural response to the problem of many non-extensions. For, if the extension operator, applied to a concept, is a sense of that concept, then it is natural to expect that there will be many senses which are not extensions. For, part of the explanatory power of the Fregean doctrine of sense is that any given referent can be presented in a number of different ways. Moreover, there is no reason to expect there to be any antecedently specified finite bound on the number of different ways that a referent can be presented. Hence, because a sense-selecting extension operator selects but one sense amongst many for each concept, there will inevitably be numerous objects in these models that are not extensions. Of course, this response to the problem of many non-extensions presupposes that one is thinking about the entities of type~$a^{\prime}$ as entities similar to Fregean senses in the respect that any given referent can be presented in a number of different ways. Even though Church himself was motivated by the project of axiomatizing  Fregean sense, obviously there is nothing written into the axioms of Church's intensional logic or the extensions thereof considered here which forces one to adopt this presupposition.

\section{Conclusions}\label{sec:RM:08}

The last two sections have illustrated some of costs and benefits of the predicative response to the Russell-Myhill paradox of propositions, at least when expanded by certain choice-like principles or by the representation operators. In the previous section~\S\ref{sec:RM:07}, we've developed a response to the Wehmeier problem of many non-extensions: the solution simply is that there are many non-extensions because the extension of a concept selects one sense from the many which present a concept. But in  section~\S\ref{sec:RM:06}, we saw that the consistency proof for the predicative response  has some features which are not in the spirit of a fine-grained theory of intensions. For, on this model, intensions of functions are only as injective as the functions they present. Perhaps there are other model constructions which would not be committed to this result, but that question is left unresolved by the work in this paper. 

Likewise, for reasons of space we have been unable to compare and contrast the versions of Church's intensional logic studied here to other formalizations of Fregean sense given by authors such as Chalmers, Horty, Moschovakis, and Tich\'y (\cite{Chalmers2011aa}, \cite{Horty2007aa}, \cite{Moschovakis1993aa}, \cite{Tichy1988aa}), or to other formal systems of fine-grained intensions due to Fox, Lappin, Parsons, and Thomason (\cite{Fox2005aa}, \cite{Parsons1982aa}, \cite{Thomason1980aa}). Finally, while we indicated in \S\ref{sec:RM:03} how others like  Anderson and Kaplan produced models which yield responses to the formalized version of the Russell-Myhill paradox~(\ref{eqn:RM:formalized}), we do not pretend to have done any serious appraisal of the costs and benefits of these proposals as compared with the predicative response. Rather, we have limited ourselves here to merely setting out the predicative response in a clear manner. In addition to suggesting a motivation for the restriction on the comprehension schema, our efforts in this paper have been directed towards establishing the formal consistency of the predicative response to the Russell-Myhill paradox of propositions. 

\section{Appendix 1: Proof of Proposition on Location of Domains}\label{app1}

In this brief appendix, we prove the Proposition on Location of Domains~(\ref{prop:location:domain}) from \S\ref{sec:RM:05}. For ease of reference, we restate it here:
\begin{myenumerate}
\item[(\ref{prop:location:domain})] For all~$n\geq 1$, both of the following hold:
\begin{itemize}
\item[] (I) for all types~$a$ with~$\|a\|<n$, there is a~$\Sigma_1$-formula in parameter~$\mu_n$ such that~$D_a$ is the unique element of~$L_{\alpha_n}$ which satisfies this formula, wherein~$\mu_{n}$ is defined by~$\mu_{n}= \langle \nu_1, \ldots, \nu_{n}, \alpha_0, \alpha_1, \ldots, \alpha_{n-1}\rangle$.
\item[] (II) for all types~$a$ with~$\|a\|= n$, the set~$D_a$ is a~$\utilde{\Sigma}_{\ell_{n}}$-definable subset of~$L_{\alpha_{n}}$ in parameter~$\mu_{n}$.
\end{itemize}
\end{myenumerate}
The proof is by simultaneous induction on~$n\geq 1$. For~$n=1$, note that (I) holds vacuously. As for (II), first note that if~$a$ is a type with~$\|a\|=1$, then~$a$ is among the types~$e,e^{\prime}, e^{\prime\prime}, \ldots, t, t^{\prime},t^{\prime\prime}, \ldots$. Now, if~$a=e$ or~$a=t$, then part~(II) follows trivially since~$\alpha_0<\alpha_1$ and so both~$\alpha_0$ and the set~$\{0,1\}$ are members of~$L_{\alpha_1}$ and~$\mu_1$ includes the parameter~$\alpha_0$ by definition. Suppose the result holds for~$a$. Since~$\mathcal{O}_{1}$ is a~$\utilde{\Sigma}_{\ell_{1}}$-definable subset of~$L_{\alpha_{1}}$ in parameter~$\nu_{1}$, it follows trivially that~$D_{a^{\prime}}=\mathcal{O}_{1}$ is~$\utilde{\Sigma}_{\ell_{1}}$-definable subset of~$L_{\alpha_{1}}$ in the more complex parameter~$\mu_{1}$. This completes the argument in the case~$n=1$.

Now suppose that the result holds for~$n$, and we show it holds for~$n\mbox{+}1$. For~(I), suppose that~$a$ is a type with~$\|a\|<n\mbox{+}1$, say~$\|a\|=m$. Then since~$D_a$ is a definable subset of~$L_{\alpha_m}$ by the induction hypothesis on part~(II) for~$m$, we may write~$D_a = \{x\in L_{\alpha_m}: L_{\alpha_m}\models \psi(x,\mu_m)\}$ for some formula~$\psi$. Hence $D_a$ is an element of~$L_{\alpha_m+1}$ and a member of~$L_{\alpha_{n+1}}$. Then we have that~$D_a$ is the unique~$X$ in~$L_{\alpha_{n+1}}$ which satisfies the following condition:
\begin{myequation}
(\forall \; x\in X\cap L_{\alpha_m} \;  L_{\alpha_m}\models \psi(x,\mu_m)) \; \& \; (\forall \; x\in L_{\alpha_m} \; (L_{\alpha_m}\models \psi(x,\mu_m)\rightarrow x\in X))
\end{myequation}
Then since~$m<n\mbox{+}1$ and the parameter~$\mu_{n+1}$ contains the parameter~$\mu_m$ as well as the ordinal~$\alpha_m$, and since the map~$\beta\mapsto L_{\beta}$ is~$\Delta_1$ in~$L_{\alpha_{n+1}}$ (cf. \cite{Devlin1984aa} II.2.8 p. 70) and since the satisfaction relation is likewise~$\Delta_1$ (cf. \cite{Devlin1984aa} I.9.10 p. 41), this is a~$\Sigma_1$-condition in~$L_{\alpha_{n+1}}$ in parameter~$\mu_{n+1}$. Here we're also appealing tacitly to the fact that the~$\Sigma_1$-formulas are closed under bounded quantification in models~$L_{\alpha}$ which satisfy~$\Sigma_1$-collection (cf. \cite{Devlin1984aa} Lemma I.11.6 p. 53). This completes the induction step for part~(I) of of the proposition. 

For the induction step for part~(II), note that the types with degree~$n\mbox{+}1$ are of the form~$a^{\prime}$ or~$ab$. Then we may do a subinduction on complexity of type. First suppose that~$\|a^{\prime}\|=n\mbox{+}1$ and suppose that the result holds for~$a$; we show it holds for~$a^{\prime}$. Since~$\mathcal{O}_{n+1}$ is~$\utilde{\Sigma}_{\ell_{n+1}}$-definable subset of~$L_{\alpha_{n+1}}$ in parameter~$\nu_{n+1}$, it follows trivially that~$D_{a^{\prime}}=\mathcal{O}_{n+1}$ is a~$\utilde{\Sigma}_{\ell_{n+1}}$-definable subset of~$L_{\alpha_{n+1}}$ in the more complex parameter~$\mu_{n+1}$. 

Second suppose that~$\|ab\|=n\mbox{+}1$, and suppose that the result holds for~$a,b$; we show it holds for~$ab$. There are two subcases here. In the first subcase, suppose that~$\|a\|\geq \|b\|$. Then by the definition of degree in~(\ref{eqn:RM:degtypinitial}), we have that~$\|a\|,\|b\|<\|ab\|$. Then if we let~$\mathrm{fnct}(f)$ abbreviate the~$\Sigma_0$-formula expressive of the graph~$f$ being functional, and fixing similar~$\Sigma_0$-definitions of~$\mathrm{dom}(f)=X$ and~$\mathrm{rng}(f)\subseteq Y$, then the set~$({D_b}^{D_a}) \cap L_{\alpha_{n+1}}$ is equal to
\begin{myequation}
\{f\in L_{\alpha_{n+1}}: \mathrm{fnct}(f) \; \& \; \exists \; X,Y \;  X=D_a \; \& \; \mathrm{dom}(f)=X \; \& \; Y = D_b \; \& \; \mathrm{rng}(f)\subseteq Y  \}
\end{myequation}
Then by part~(I), we have that this is a~$\utilde{\Sigma}_1$-definable subset of~$L_{\alpha_{n+1}}$ in parameter~$\mu_{n+1}$. 

Now, as a second subcase, suppose the result holds for~$a,b$ and that~$\|a\|< \|b\|$, so that by the definition of degree in~(\ref{eqn:RM:degtypinitial}) we have~$n\mbox{+}1=\|ab\|=\|b\|$. Then~$D_a$ is a member of~$L_{\alpha_{n+1}}$ by part~(I), while by the supposition that the result holds for~$b$ we have that~$D_b$ is a~$\utilde{\Sigma}_{\ell_{n+1}}$-definable subset of~$L_{\alpha_{n+1}}$ in parameter~$\mu_{n+1}$. Then~$({D_b}^{D_a}) \cap L_{\alpha_{n+1}}$ is also a~$\utilde{\Sigma}_{\ell_{n+1}}$-definable subset of~$L_{\alpha_{n+1}}$ in the parameter~$\mu_{n+1}$. For, we have the following definition of~$({D_b}^{D_a}) \cap L_{\alpha_{n+1}}$:
\begin{myequation}
\{f\in L_{\alpha_{n+1}}: L_{\alpha_{n+1}}\models [\mathrm{fnct}(f) \; \& \; \exists \; X \; X=D_a \; \& \; \forall \; x\in X \; \exists \; y\in D_b \; \langle x,y\rangle\in f]\}
\end{myequation}
Here we are appealing to part~(I) applied to~$D_a$ since~$\|a\|<n\mbox{+}1$ in this subcase. Likewise, we are appealing to the fact that the bounded quantification in the last conjunct does not move us out of the complexity class~$\utilde{\Sigma}_{\ell_{n+1}}$ in models of~$\utilde{\Sigma}_{\ell_{n+1}}$-collection and~$\utilde{\Sigma}_{\ell_{n+1}-1}$-separation. This finishes the induction step for~(II). With this the inductive proof of the proposition is finished. 
 \section{Appendix 2: Verification of the Satisfaction of Predicative Comprehension}\label{app2}

Here we prove the following theorem from \S\ref{sec:RM:05}:
\begin{myenumerate}
\item[(\ref{thM:RM:big})] (Theorem on Consistency of Predicative Comprehension) For every intensional hierarchy~$D$~(\ref{eqn:RM:defn:intensional:hierarchy}), the associated intensional structure~$\mathbb{D}$ models each instance of the Predicative Typed Choice Schema~(\ref{eqn:RM:compschem4}) and hence each instance of the Predicative Typed Comprehension Schema~(\ref{eqn:RM:predcompschem}).
\end{myenumerate}
Further, we here prove this result for the language expanded by the representation functions~$\nabla_a: D_a\rightarrow D_{a^{\prime}}$ introduced in \S\ref{sec:RM:06} (cf. circa~(\ref{eqn:RM:PresentationRepresentation})). As a first step towards approaching the proof of this theorem, let's first note an elementary result on terms. The terms in the signature of an intensional structures consists simply of the closure of the constants~$0$,$1$ and the variables under the extensional application symbols and the representation operations. The presentation symbols and the intensional application symbols are not total and hence are formally treated as relation symbols as opposed to function symbols. The \emph{type} of a term is defined inductively as follows: the truth-values~$0$,~$1$ have type~$t$, the variables have the type that they are given initially, and if~$\tau$ has type~$ab$ and~$\sigma$ has type~$a$, then~$\tau(\sigma)$ or~$\mathrm{e\mbox{-}app}_{ab}(\tau, \sigma)$ has type~$b$; and if~$\tau$ has type~$a$ then~$\nabla_a(\tau)$ has type~$a^{\prime}$. Then we have the following elementary result:
\begin{myenumerate}
\item (Proposition that Terms do not Raise Degree). Suppose that~$\tau(x_1, \ldots, x_k)$ is a term in the signature of intensional structures with all free variables displayed such that the type of each variable~$x_i$ has degree~$\leq~n$. Then the type of the term~$\tau$ has degree~$\leq~n$. \label{eqn:RM:dontraise}\index{Proposition that Terms do not Raise Degree (\ref{eqn:RM:dontraise})}
\end{myenumerate}
The proof is by induction on the complexity of the term. Clearly this is true in the case of the truth-values and the variables. Suppose it holds for~$\tau(x_1, \ldots, x_k)$ and~$\sigma(x_1, \ldots, x_k)$; we must show it is the case for~$\mathrm{e\mbox{-}app}_{ab}(\tau, \sigma)$ and~$\nabla_a(\tau)$. First consider the case of~$\mathrm{e\mbox{-}app}_{ab}(\tau, \sigma)$. Then~$\tau$ has type~$ab$ and~$\sigma$ has type~$a$, and each has type with degree~$\leq n$ by the induction hypothesis. There are two cases to consider, corresponding to the two clauses in the definition of~$\|ab\|$ in~(\ref{eqn:RM:degtypinitial}). First suppose that~$\|a\|\geq \|b\|$. Then one has that~$\|b\|\leq \|a\|\leq n$, which is what we wanted to show since the type of~$\mathrm{e\mbox{-}app}_{ab}(\tau, \sigma)$ is~$b$. Second suppose that~$\|a\|<\|b\|$. Then we have that~$\|b\|=\|ab\|\leq n$, which is again what we wanted to show. Finally, consider the case of~$\nabla_a(\tau)$. Then~$\tau$ has type~$a$, and it has degree~$\leq{n}$ by induction hypothesis. Then~$\nabla_a(\tau)$ has type~$a^{\prime}$ and so~$\|a^{\prime}\|=\|a\|\leq n$ by the definition of degree of~$\|a^{\prime}\|$ in~(\ref{eqn:RM:degtypinitial}). This is why terms do not raise degree, or why~(\ref{eqn:RM:dontraise}) holds.

Relatedly, as a preliminary step, let's establish the following result about the complexity of the functions on intensional structures induced by terms:
\begin{myenumerate}
\item (Proposition on Complexity of Terms) Suppose that~$\tau(\overline{u})\equiv \tau(u_1, \ldots, u_j)$ is a term with all free variables displayed where~$u_i$ has type~$d_i$. Since terms don't raise degree~(\ref{eqn:RM:dontraise}),~$\tau$ has type with degree~$d$ with~$\|d\|\leq m=\max\{\|d_1\|, \ldots, \|d_j\|\}$. Then~$\tau$ induces a function~$\tau^{\mathbb{D}}: D_{d_1}\times \cdots \times D_{d_j}\rightarrow D_d$ whose graph is~$\utilde{\Sigma}^{L_{\alpha_m}}_{\ell_m}$-definable.\label{eqn:asdfasdfdsaafasdf3214231}\index{Proposition on Complexity of Terms (\ref{eqn:asdfasdfdsaafasdf3214231})}
\end{myenumerate} 
Clearly this is the case if the term is variable. Now for the induction step suppose that the result holds for~$\tau$ and~$\sigma$; we must show it holds for~$\rho(\overline{u})\equiv \mathrm{e\mbox{-}app}_{e_1e_2}(\tau(\overline{u}), \sigma(\overline{u}))$. Then~$\tau(\overline{u})$ has type~$e_1e_2$ and~$\sigma(\overline{u})$ has type~$e_1$. Since terms don't raise degree~~(\ref{eqn:RM:dontraise}), it follows that~$\|e_1e_2\|, \|e_1\|$ are all less than or equal to~$m=\max\{\|d_1\|, \ldots, \|d_j\|\}$, and from this we infer that~$\|e_2\|\leq \|e_1e_2\|\leq m$ as well. Then~$\rho^{\mathbb{D}}: D_{d_1}\times \cdots \times D_{d_j}\rightarrow D_{e_2}$ has the following graph:
\begin{myeqnarray}
\{ (\overline{u},u)\in D_{d_1}\times \cdots \times D_{d_j}\times D_{e_2} & : &  \exists \; y\in D_{e_1} \; \exists \; z \in D_{e_1e_2} \;  \notag \\
& & \sigma^{\mathbb{D}}(\overline{u})= y\; \& \;  \tau^{\mathbb{D}}(\overline{u})=z \; \& \; \langle y,u\rangle \in z\} 
\end{myeqnarray}
This is~$\utilde{\Sigma}^{L_{\alpha_m}}_{\ell_m}$-definable by the Location of Domains~(\ref{prop:location:domain}) since we have that~$\|e_1\|, \|e_2\|,\|e_1e_2\|\leq m$.  For the final induction step, suppose that the result holds for~$\tau$; we must show it holds for~$\nabla_a(\tau)$. Then~$\tau$ has type~$a$, and since terms don't raise degree~~(\ref{eqn:RM:dontraise}), it follows that~$\|a\|\leq m$. Then the graph of~$\tau^{\mathbb{D}}$ is~$\utilde{\Sigma}^{L_{\alpha_m}}_m$-definable by induction hypothesis. Recall from the discussion in \S\ref{sec:RM:06} that the representation function~$\nabla_a$ is interpreted on intensional structures~$\mathbb{D}$ by the function~$\iota_{\|a\|}$ from the definition of an intensional hierarchy~(\ref{eqn:RM:defn:intensional:hierarchy}). However, this was by definition~$\utilde{\Sigma}^{L_{\alpha_{\|a\|}}}_{\|a\|}$-definable (cf. clause (iii) of the definition of an intensional position~(\ref{eqn:RM:defn:intensional:position})). Since~$\|a\|\leq m$, we then have that the composition~$\iota_{\|a\|}\circ \tau^{\mathbb{D}}$ is clearly also ~$\utilde{\Sigma}^{L_{\alpha_m}}_m$-definable. This finishes the proof of result on the complexity of terms~(\ref{eqn:asdfasdfdsaafasdf3214231}).

Now let's consider what kinds of symbols can appear in a formula covered by the Predicative Typed Choice Schema~(\ref{eqn:RM:compschem4}). Suppose that the formula~$\varphi(x, y, z_1, \ldots, z_k)$ is a formula with all free variables displayed and with free variable~$x$ of type~$a$,~$y$ of type~$b$, and in addition variable~$z_i$ has type~$c_i$ with~$\|c_i\|\leq \|ab\|$ and all the bound variables in~$\varphi(x,y,z_1, \ldots, z_k)$ have type~$c$ with~$\|c\|< \|ab\|$. Let~$\|ab\|=n+1$. There are then two cases to consider, corresponding to the split in cases in the definition of the degree~$\|ab\|$ in~(\ref{eqn:RM:degtypinitial}). If~$\|a\|\geq \|b\|$, then~$n+1=\|ab\|=\|a\|+1$ and so~$\|b\|\leq \|a\|\leq n$. Further, if we split the parameter variables~$z_1, \ldots, z_k$ into those that have type with degree~$n+1$ and those that have type with degree~$\leq~n$, then we can write the formula in question as:
\begin{myenumerate}
\item (First Configuration):~$\varphi(x, y, v_1, \ldots, v_m, z_1, \ldots, z_k)$ is a formula with all free variables displayed and with free variable~$x$ of type~$a$ with~$\|a\|\leq{n}$,~$y$ of type~$b$ with~$\|b\|\leq{n}$, and in addition variable~$v_i$ has type~$a_i$ with~$\|a_i\|\leq{n}$ and variable~$z_i$ has type~$c_i$ with~$\|c_i\|=n+1$ and all the bound variables in the formula have type~$c$ with~$\|c\|\leq{n}$.\label{eqn:RM:firstconfig}\index{Configuration, First (\ref{eqn:RM:firstconfig})}
\end{myenumerate}
Alternatively, in the other case, we have~$\|a\|<\|b\|$ and~$n+1=\|ab\|=\|b\|$. If we again split the parameter variables~$z_1, \ldots, z_k$ into those that have type with degree~$n+1$ and those that have type with degree~$\leq~n$, then we can write the formula in question as:
\begin{myenumerate}
\item (Second Configuration):~$\varphi(x, v_1, \ldots, v_m, y, z_1, \ldots, z_k)$ is a formula with all free variables displayed and with free variable~$x$ of type~$a$ with~$\|a\|\leq{n}$,~$y$ of type~$b$ with~$\|b\|= n+1$, and in addition variable~$v_i$ has type~$a_i$ with~$\|a_i\|\leq{n}$ and variable~$z_i$ has type~$c_i$ with~$\|c_i\|=n+1$ and all the bound variables in the formula have type~$c$ with~$\|c\|\leq{n}$.\label{eqn:RM:secondconfig}\index{Configuration, Second (\ref{eqn:RM:secondconfig})}
\end{myenumerate}
For ease of future reference, we call these two kinds of formulas which can feature in the Predicative Typed Choice Schema~(\ref{eqn:RM:compschem4}) the ``first configuration'' and the ``second configuration''.

The plan in what follows is to show that the Predicative Typed Choice Schema~(\ref{eqn:RM:compschem4}) holds for formulas in the second configuration~(\ref{eqn:RM:secondconfig}), and then to show it for formulas in the first configuration~(\ref{eqn:RM:firstconfig}). This first step is done by proving a result connecting the satisfaction of a formula in the second configuration to a certain level of definability in the constructible hierarchy. To build up to the statement of this result, suppose that~$\varphi(x, v_1, \ldots, v_m, y, z_1, \ldots, z_k)$ is in the second configuration~(\ref{eqn:RM:secondconfig}). Then any subformula of this formula has the form
\begin{myequation}
\psi(x,v_1, \ldots, v_m, v_{m+1}, \ldots v_{m+m^{\prime}},y,z_1, \ldots, z_k)
\end{myequation}
where the variable~$v_i$ for~$i>m$ has type~$a_i$ with degree~$\leq{n}$. Let's abbreviate~$\overline{v} = \langle v_1, \ldots, v_m, v_{m+1}, \ldots v_{m+m^{\prime}}\rangle$ and let's abbreviate 
\begin{myequation}
D_{\overline{a}}=D_{a_1}\times \cdots \times D_{a_{m+m^{\prime}}}, \hspace{10mm} D_{\overline{c}} = D_{c_1}\times \cdots \times D_{c_k}
\end{myequation}
Note that since~$\|a\|, \|a_i\|\leq n$, it follows from the Location of Domains~(\ref{prop:location:domain}), we have that~$D_{a}\times D_{\overline{a}}$ is a member of~$L_{\alpha_{n+1}}$. However, since~$\|b\|, \|c_i\|=n+1$, we have that~$D_b\times D_{\overline{c}}$ is a~$\utilde{\Sigma}^{L_{\alpha_{n+1}}}_{\ell_{n+1}}$-definable subset of~$L_{\alpha_{n+1}}$.  Having put this terminology in place, let's now show that:
\begin{myenumerate}
\item (Proposition on Complexity of Satisfaction, Second Configuration) For every intensional hierarchy~$D$ with induced intensional structure~$\mathbb{D}$ and every subformula~$\psi(x,\overline{v},y,\overline{z})$ of a formula in the second configuration~(\ref{eqn:RM:secondconfig}), the following set is~$\utilde{\Sigma}_{\ell_{n+1}}^{L_{\alpha_{n+1}}}$-definable:\label{eqn:RM:complexsat:second}\index{Proposition on Complexity of Satisfaction, Second Configuration (\ref{eqn:RM:complexsat:second})}
\[[ \psi]^D= \{(x,\overline{v},y,\overline{z})\in D_a\times D_{\overline{a}}\times D_b\times D_{\overline{c}} :  \mathbb{D}\models \psi(x,\overline{v},y,\overline{z})\}\] \end{myenumerate}
We establish this by induction on the complexity of the subformula. By pushing all the negations to the inside, it suffices to show that the result holds for atomics, negated atomics, and is closed under conjunctions, disjunctions, existential quantification, and universal quantification. Let's begin with the atomic case, considering the negated atomic cases along the way. The atomic formulas in intensional structures have three possible forms, namely:
\begin{myequation}
\tau = \sigma, \hspace{10mm} \Delta_d(\tau)=\sigma, \hspace{10mm} \mathrm{i\mbox{-}app}_{c_0d_0}(\tau, \sigma) = \rho
\end{myequation}
where~$\tau, \sigma, \rho$ are terms. These are the only possible subformulas because technically, the second is shorthand for the binary atomic formula~$\Delta_d(\tau, \sigma)$ and the third is shorthand for the associated ternary atomic relation (cf. discussion circa equations~(\ref{eqn:RM:defnpres}) and (\ref{eqn:RM:defn:intensional})). Since~$\tau, \sigma,\rho$ appear in a formula in the second configuration~(\ref{eqn:RM:secondconfig}), the free variables in these terms~$\tau, \sigma, \rho$ have types with degree~$\leq{n}+1$ and since terms don't raise degree~(\ref{eqn:RM:dontraise}), it follows that the respective types~$e_1, e_2, e_3$ of~$\tau, \sigma, \rho$ are also such that~$\|e_1\|, \|e_2\|, \|e_3\| \leq{n}+1$. From this it follows in turn that~$\|d\|\leq{n}+1$ and~$\|c_0\|, \|d_0\| \leq \|c_0d_0\|=\|(c_0d_0)^{\prime}\|\leq{n}+1$. 

Let's consider first the case of equality between terms, that is, atomic formulas of the form
\begin{myequation}
\psi(x,\overline{v}, y,\overline{z})\equiv \tau(x,\overline{v},y,\overline{z})=\sigma(x,\overline{v},y, \overline{z})
\end{myequation}
Then we have that
\begin{myeqnarray}
(x,\overline{v},y,\overline{z})\in [\psi]^D & \Longleftrightarrow &  \exists \; z_1\in D_{e_1}, \exists \; z_2\in D_{e_2} \; \tau(x,\overline{v},y,\overline{z})=z_1\notag \\
&& \; \& \; \sigma(x,\overline{v},y,\overline{z})=z_2 \; \& \; z_1=z_2 
\end{myeqnarray}
which is~$\utilde{\Sigma}_{\ell_{n+1}}^{L_{\alpha_{n+1}}}$-definable by the result on the complexity of terms~(\ref{eqn:asdfasdfdsaafasdf3214231}). Similarly we have that
\begin{myeqnarray}
(x,\overline{v},y, \overline{z})\in [\neg \psi]^D & \Longleftrightarrow &  \exists \; z_1\in D_{e_1}, \exists \; z_2\in D_{e_2} \; \tau(x,\overline{v},y, \overline{z})=z_1 \notag \\
& &  \; \& \; \sigma(x,\overline{v},y, \overline{z})=z_2 \; \& \; z_1\neq z_2
\end{myeqnarray}
which is~$\utilde{\Sigma}_{\ell_{n+1}}^{L_{\alpha_{n+1}}}$-definable for the same reasons. 

Now let's consider the case of the presentation symbols, that is atomic formulas of the form
\begin{myequation}
\psi(x,\overline{v},y, \overline{z})\equiv \Delta_d(\tau(x,\overline{v},y,\overline{z}))=\sigma(x,\overline{v},y,\overline{z})
\end{myequation}
Then~$[\psi]^D$ is~$\utilde{\Sigma}_{\ell_{n+1}}^{L_{\alpha_{n+1}}}$-definable because we have the following biconditional and because~$\|d\|\leq{n}\mbox{+}1$ implies that~$\pi_{\|d\|}$ is~$\utilde{\Sigma}_{\ell_{n+1}}^{L_{\alpha_{n+1}}}$-definable:
\begin{myeqnarray}
(x,\overline{v},y,\overline{z})\in [\psi]^D & \Longleftrightarrow &  \exists \; z_1\in D_{e_1}, \exists \; z_2\in D_{e_2} \; \tau(x,\overline{v},y,\overline{z})=z_1\notag \\
& & \; \& \; \sigma(x,\overline{v},y,\overline{z})=z_2 \; \& \; \pi_{\|d\|}(z_1)=z_2 
\end{myeqnarray}
Further, by using the fact that $\mathcal{O}_j\setminus \pi^{-1}_j(L_{\alpha_j})$ was~$\utilde{\Sigma}_{\ell_{j}}^{L_{\alpha_{j}}}$-definable for all~$j\geq 1$  (cf. clause~(v) of the definition of an intensional position~(\ref{eqn:RM:defn:intensional:position})) and that~$\|d\|\leq{n}+1$, we have that 
\begin{myeqnarray}\label{eqn:d123412341234213}
(x,\overline{v},y,\overline{z})\in [\neg\psi]^D & \Longleftrightarrow&  \exists \; z_1\in D_{e_1}, \exists \; z_2 \in D_{e_2}, \; \exists \; z_3\in L_{\alpha_{n+1}}\notag \\
& & \tau(x,\overline{v},y,\overline{z})=z_1  \; \& \; \sigma(x,\overline{v},y, \overline{z})=z_2 \\
& &  \; \& \; (z_1\in \mathcal{O}_{\|d\|}\setminus \pi^{-1}_{\|d\|}(L_{\alpha_{\|d\|}})) \vee (\pi_{\|d\|}(z_1)=z_3 \; \& \; z_2\neq z_3)\notag
\end{myeqnarray}

As the final atomic case, consider the case of intensional application:
\begin{myequation}
\psi(x,\overline{v},y, \overline{z})\equiv \mathrm{i\mbox{-}app}_{c_0d_0}(\tau(x,\overline{v},y,\overline{z}), \sigma(x,\overline{v},y,\overline{z})) = \rho(x,\overline{v},y,\overline{z})
\end{myequation}
Then the type~$e_1$ of~$\tau$ must be~$(c_0d_0)^{\prime}$ and the type~$e_2$ of~$\sigma$ must be~$c_0^{\prime}$. Then we have
\begin{myeqnarray}
(x,\overline{v},y,\overline{z})\in [\psi]^D & \Longleftrightarrow & \exists \; z_1\in D_{(c_0d_0)^{\prime}}, \exists \; z_2\in D_{c_0^{\prime}}, \exists \; z_3\in D_{d_0^{\prime}}\notag \\
& &  \tau(x,\overline{v},y,\overline{z}) = z_1 \; \& \;  \sigma(x,\overline{v},y,\overline{z})=z_2 \; \& \; \rho(x,\overline{v},y,\overline{z})=z_3 \notag \\
& & \; \& \; \iota_{\|d_0\|}((\pi_{\|c_0d_0\|}(z_1))(\pi_{\|c_0\|}(z_2)) = z_3 
\end{myeqnarray}
Since~$\|c_0\|, \|d_0\| \leq \|c_0d_0\|=\|(c_0d_0)^{\prime}\|\leq{n}+1$, we have that~$\pi_{\|c_0d_0\|}, \pi_{\|c_0\|}, \iota_{\|d_0\|}$ are all~$\utilde{\Sigma}_{\ell_{n+1}}^{L_{\alpha_{n+1}}}$-definable. Finally, for the negation, we may argue as follows, again appealing to the fact that in the definition of an intensional position we required that the set~$\mathcal{O}_j\setminus \pi^{-1}_j(L_{\alpha_j})$ was~$\utilde{\Sigma}_{\ell_{j}}^{L_{\alpha_{j}}}$-definable for all~$j\geq 1$  (cf. clause~(v) of the definition of an intensional position~(\ref{eqn:RM:defn:intensional:position})):
\begin{myeqnarray}\label{eqn:complemee321423}
(x,\overline{v},y,\overline{z})\in [\neg \psi]^D & \Longleftrightarrow & \exists \; z_1\in D_{(c_0d_0)^{\prime}}, \exists \; z_2 \in D_{c_0^{\prime}}, \exists \; z_3\in D_{d_0^{\prime}}, \exists \; z_4\in L_{\alpha_{n+1}} \notag \\
& &  \tau(x,\overline{v},y,\overline{z}) = z_1 \; \& \;  \sigma(x,\overline{v},y,\overline{z})=z_2 \; \& \; \rho(x,\overline{v},y,\overline{z})=z_3 \notag \\
& & \wedge [(z_1\in \mathcal{O}_{\|c_0d_0\|}\setminus \pi^{-1}_{\|c_0d0\|}(L_{\alpha_{\|c_0d_0\|}}))\notag \\  &&\;  \vee\;  (z_2\in \mathcal{O}_{\|c_0\|}\setminus \pi^{-1}_{\|c_0\|}(L_{\alpha_{\|c_0\|}}))    \\
& & \; \vee \; (\iota_{\|d_0\|}((\pi_{\|c_0d_0\|}(z_1))(\pi_{\|c_0\|}(z_2)) = z_4 \; \& \; z_4\neq z_3)] \notag
\end{myeqnarray}

This completes the base cases of the inductive argument for~(\ref{eqn:RM:complexsat:second}). Since~$\utilde{\Sigma}_{\ell_{n+1}}^{L_{\alpha_{n+1}}}$-definability is closed under finite intersections and unions, the inductive steps for conjunction and disjunction are trivial. Let us then consider the case of universal quantification. Suppose that the result holds for~$\psi(x,\overline{v},v,y,\overline{z})$ and let us show it holds for~$\theta(x,\overline{v},y,\overline{z})\equiv \forall \; v_0 \; \psi(x,\overline{v},v_0,y,\overline{z})$. Since this is a subformula of a formula in the second configuration~(\ref{eqn:RM:secondconfig}), it follows the bound variable~$v_0$ has a type~$a_0$ with degree~$\leq{n}$. Then by part~(ii) of the result on Locations of Domains~(\ref{prop:location:domain}), it follows that~$X=D_{a_0}$ is a~$\utilde{\Sigma}^{L_{\alpha_{n+1}}}_1$-condition. Then one has that
\begin{myequation}
(x,\overline{v},y,\overline{z})\in [\theta]^D \Longleftrightarrow \exists \; X \; X = D_{a_0} \; \& \; \forall \; v_0\in X \; (x,v_0,\overline{v}, y,\overline{z})\in [\psi]^D
\end{myequation}
so that~$[\theta]^D$ is likewise ~$\utilde{\Sigma}_{\ell_{n+1}}^{L_{\alpha_{n+1}}}$-definable since ~$\utilde{\Sigma}_{\ell_{n+1}}^{L_{\alpha_{n+1}}}$-definability is closed under bounded quantification in models~$L_{\alpha_{n+1}}$ of~$\Sigma_{\ell_{n+1}}$-collection and~$\Sigma_{\ell_{n+1}-1}$-separation. A similar argument holds in the case of the existential quantifier, but is even easier since there we do not have to appeal to this result about closure under bounded quantification. This finishes the result on the complexity of satisfaction in the case of a formula which is in the second configuration~(\ref{eqn:RM:complexsat:second}).

Now let us finally establish that the Predicative Typed Choice Schema~(\ref{eqn:RM:predcompschem}) holds on intensional structures, at first with respect to formulas in the second configuration~(\ref{eqn:RM:secondconfig}). Suppose that the antecedent holds:
\begin{myequation}\label{eqn:testmealready}
\mathbb{D}\models \forall \;x \; \exists \; y \; \varphi(x,p_1, \ldots, p_m, y,q_1, \ldots, q_k)
\end{myequation}
where~$\varphi(x,v_1, \ldots,v_m, y,z_1, \ldots z_k)$ is in the second configuration~(\ref{eqn:RM:secondconfig}). Consider the following relation:
\begin{myequation}
R(x,y)\equiv [x\in D_a \; \& \; y\in D_b \; \& \; \mathbb{D}\models \varphi(x,p_1, \ldots, p_m,y, q_1, \ldots, q_k)]
\end{myequation}
Then by the result on the complexity of satisfaction~(\ref{eqn:RM:complexsat:second}), one has that~$R$ is~$\utilde{\Sigma}_{\ell_{n+1}}^{L_{\alpha_{n+1}}}$-definable. And by equation~(\ref{eqn:testmealready}), one has that 
\begin{myequation}\label{eqn:RM:almostdone}
L_{\alpha_{n+1}}\models \forall \;x\in D_a \; \exists \; y \; R(x,y)
\end{myequation}
By the uniformization theorem (cf. \cite{Jensen1972aa} Theorem 3.1 p. 256 and Lemma 2.15 p. 255; \cite{Devlin1984aa} Theorem 4.5 p. 269, and ``weak uniformization'' in \cite{Walsh2014ac}), choose a~$\utilde{\Sigma}_{\ell_{n+1}}^{L_{\alpha_{n+1}}}$-definable relation~$R^{\prime}\subseteq R$ such that 
\begin{myequation}
L_{\alpha_{n+1}} \models [\forall \; x \; (\exists \; y \; R(x,y))\rightarrow (\exists \; ! \; y \; R^{\prime}(x,y))]
\end{myequation}
Then by equation~(\ref{eqn:RM:almostdone}), one has that~$R^{\prime}$ is the graph of a function~$h:D_a\rightarrow D_b$. Since this graph is ~$\utilde{\Sigma}_{\ell_{n+1}}^{L_{\alpha_{n+1}}}$-definable with domain~$D_a$ an element of~$L_{\alpha_{n+1}}$, by Replacement (cf. \cite{Devlin1984aa} Lemma I.11.7 p. 53) one has that it is an element of~$L_{\alpha_{n+1}}=L_{\alpha_{\|ab\|}}$. Then~$h$ is an element of the domain~$D_b^{D_a}\cap L_{\alpha_{\|ab\|}}= D_{ab}$ (cf. the third clause of equation~(\ref{eqn:defn:RM:typestodomains})). Hence, we've shown that there is~$h$ in~$D_{ab}$ such that 
\begin{myequation}
\mathbb{D}\models \forall \;x \; \varphi(x,p_1, \ldots, p_m, h(x),q_1, \ldots, q_k)
\end{myequation}
which is what we were required to show in the consequent of the Predicative Typed Choice Schema~(\ref{eqn:RM:predcompschem}).

\vspace{2mm}

We've verified that the Predicative Typed Choice Schema~(\ref{eqn:RM:predcompschem}) holds on intensional structures, at least with respect to formulas in the second configuration~(\ref{eqn:RM:secondconfig}). Let's now argue that the same holds with respect to formulas in the first configuration~(\ref{eqn:RM:firstconfig}). Suppose that
\begin{myequation}\label{eqn:testmealready2}
\mathbb{D}\models \forall \;x \; \exists \; y \; \varphi(x,y,p_1, \ldots, p_m, q_1, \ldots, q_k)
\end{myequation}
where~$\varphi(x,y,v_1, \ldots,v_m, z_1, \ldots z_k)$ is in the first configuration~(\ref{eqn:RM:firstconfig}). Then consider the following, where~$w$ is a variable of type~$ab$ with degree~$n+1$ and~$x_1,x_2$ are variables of type~$a$:
\begin{myequation}
\psi(x,w,\overline{p}, \overline{q})\equiv (\forall \; x_1, x_2\; w(x_1)=w(x_2)) \; \& \; (\exists \; x_1, y \;  \; (w(x_1)=y \; \& \; \varphi(x,y,\overline{p}, \overline{q})))
\end{myequation}
Intuitively this formula~$\psi$ is saying that~$w$ is a constant function of type~$ab$ and its constant value is a witness to~$\varphi$. Now~$\psi$ is in the second configuration~(\ref{eqn:RM:secondconfig}), and we can verify by hand that for every element~$y$ of~$D_b$ there is a constant function of type~$ab$ whose constant value is $y$. For, if~$y\in L_{\alpha_{\|b\|}}$ then~$\{\langle x_1,y\rangle: x_1\in D_a\}$ is in~$L_{\alpha_{\|ab\|}}$. Then by  Predicative Typed Choice Schema~(\ref{eqn:RM:predcompschem}) applied to~$\psi$, we have that there is an element~$h$ of type~$a(ab)$ such that~$\mathbb{D}\models \forall \; x \; \psi(x,h(x), p_1, \ldots, p_m, q_1, \ldots, q_k)$. Note that since~$\|ab\|=n+1=\|a(ab)\|$ we have that~$h$ is in~$L_{\alpha_{n+1}}$. Then the function~$g:D_a\rightarrow D_b$ such that~$g(x)=(h(x))(x)$ is in~$L_{\alpha_{n+1}}$ by~$\Sigma_0$-separation since  
\begin{myequation}
g =  \{\langle x,y\rangle\in D_a\times D_b: \langle x, \langle x,y\rangle\rangle \in h\}
\end{myequation}
We've now shown that~$g$ is an element of~$D_{ab}$ and by construction we have  
\begin{myequation}\label{eqn:testmealready3}
\mathbb{D}\models \forall \;x \; \varphi(x,g(x),p_1, \ldots, p_m, q_1, \ldots, q_k)
\end{myequation}
so that we also have that Predicative Typed Choice Schema~(\ref{eqn:RM:predcompschem}) holds on intensional structures, regardless of which of the two configurations we are in.

\subsection*{Acknowledgements}
I was lucky enough to be able to present parts of this work at a number of workshops and conferences, and I would like to thank the participants and organizers of these events for these opportunities. I would like to especially thank the following people for the comments and feedback which I received on these and other occasions: Robert~Black, Roy~Cook, Matthew~Davidson, Walter~Dean, Marie~Du\v{z}\'{\i}, Kenny~Easwaran, Fernando~Ferreira, Martin~Fischer, Rohan~French, Salvatore~Florio, Kentaro~Fujimoto, Jeremy~Heis, Joel~David~Hamkins, Volker~Halbach, Ole Thomassen~Hjortland, Luca~Incurvati, Daniel~Isaacson, J\"onne~Kriener, Graham~Leach-Krouse, Hannes~Leitgeb, {\O}ystein~Linnebo, Paolo~Mancosu, Richard~Mendelsohn, Tony~Martin, Yiannis~Moschovakis, John~Mumma, Pavel~Pudl\'ak, Sam~Roberts, Marcus~Rossberg, Tony~Roy, Gil~Sagi, Florian~Steinberger, Iulian~Toader, Gabriel~Uzquiano, Albert~Visser, Kai~Wehmeier, Philip~Welch, Trevor~Wilson, and Martin~Zeman. This paper has likewise been substantially bettered by the feedback and comments of the editors and referees of this journal, to whom I express my gratitude. While composing this paper, I was supported by a Kurt G\"odel Society Research Prize Fellowship and by {\O}ystein Linnebo's European Research Council funded project ``Plurals, Predicates, and Paradox.'' 

\bibliographystyle{plain}
\bibliography{walsh-sean-final.bib}

\end{document}